\documentclass[12pt,reqno]{article}
\usepackage{amsmath,amssymb,latexsym,amsthm}
\input xy
\xyoption{all}
\CompileMatrices

\newtheorem{thm}{Theorem}[section]
\newtheorem{prop}[thm]{Proposition}
\newtheorem{lemma}[thm]{Lemma}
\newtheorem{cor}[thm]{Corollary}

\theoremstyle{definition}
\newtheorem{defn}[thm]{Definition}

\theoremstyle{remark}
\newtheorem{rem}[thm]{Remark}
\newtheorem{rems}[thm]{Remarks}

\theoremstyle{example}
\newtheorem{examples}[thm]{Examples}

\DeclareMathOperator{\tr}{Tr}
\DeclareMathOperator{\im}{Im}
\DeclareMathOperator{\id}{Id}

\let\dsp=\displaystyle

\def\R{\mathbb R}
\def\C{\mathbb C}
\def\N{\mathbb N}

\def\Q{\mathbb Q}
\def\H{\mathbb H}

\begin{document}

\title{Around heat decay on forms \\ and relations of nilpotent Lie
  groups}
\author{Michel Rumin} %\date{\today }

\date{}

\maketitle

\abstract{One knows that the large time heat decay exponent on a
  nilpotent group is given by half the growing rate of the volume of
  its large balls. This work deals with the similar problem of trying
  to interpret geometrically the heat decay on (one) forms. We will
  show how it is (partially) related to the depth of the relations
  required to define the group. The tools used apply in general on
  Carnot-Carath\'eodory manifolds.}

\tableofcontents{}

\section{Introduction}

Let $(M,g)$ be a compact riemannian manifold with fundamental group
$\Gamma =\pi_1(M)$. We denote by $\tilde M$ its universal cover and
$\tilde g$ the pull back metric. The de Rham differential
$d:\Omega^p(\tilde M) \rightarrow \Omega^{p+1} (\tilde M)$ acts
between $p$ and $p+1$-forms of $\tilde M$.

We are interested in the spectrum of the $p$-Laplacian $\Delta_p =
d\delta +\delta d$ acting on $L^2$ $p$-forms of $\tilde M$. More
precisely, one would like to know what kind of geometric information
about $M$ or $\Gamma$ is encoded in the near zero spectrum of
$\Delta_p$.

The first result of this kind deals with the invertibility of
$\Delta_0$, the Laplacian on functions.
\begin{thm}
  (Brooks \cite{Br}) $0$ belongs to the spectrum of $\Delta_0$ iff
  $\Gamma$ is amenable. (That means that $\Gamma$ may be exhausted by
  parts $\Gamma_i$ such that $\dsp
  \frac{\#(\partial\Gamma_i)}{\#(\Gamma_i)}\rightarrow 0$.)
\end{thm}
Roughly, the argument that $\Delta_0$ is not invertible if $\Gamma$ is
amenable is based on the idea that one can approach the constant
function on $\tilde M$ by compactly supported test functions with
small differential with respect to their integral norm.  These
approximative units can also be used to cut off the pull back of any
harmonic \emph{form} of $M$. This shows therefore that if $\Gamma$ is
amenable and $H^p(M,\R)\not= 0$ then $0$ belongs to the spectrum of
$\Delta_p$ on $\tilde M$. This is satisfied for any $p\leq
\mathrm{dim} M$ in the case we will study of $M=G/\Gamma$, with
$\Gamma$ a cocompact group of a nilpotent Lie group $G$ (see
\cite{Dix}).

Note also that in these examples, $0$ is never \emph{embedded} in
$\mathrm{Sp} (\Delta_p)$. There is no harmonic $L^2$-form on $\tilde
M=G$, otherwise they would have to be invariant through the Killing
direction associated to the non-vanishing center of these groups. More
precise information about the density of $\mathrm{Sp} (\Delta_p)$ near
$0$ is obtained with the help of the notions of $\Gamma$-dimension,
and $\Gamma$-trace.

Recall (see \cite{At}) that in the case of a $\Gamma$-invariant
smoothing operator $S$, acting on $\Omega^*(\tilde M)$, one can
consider the average diagonal trace density of its Schwarz kernel
$K_S$, namely
$$
\tr_\Gamma(S) = \int_{\cal F} \tr(K_S(x,x)) \mathrm{dvol}\,,$$
where ${\cal F}$ stands for any fundamental domain of the
$\Gamma$-action. For example $S$ can be the heat operator
$e^{-t\Delta_p}$ or the spectral projection $\Pi(\Delta_p \leq
\lambda)$.  Actually we can first take profit of the orthogonal
splitting of $\Delta_p$ in $\delta d$ and $d\delta$ to slightly
precise the analysis. So the basic operators in concern will be
$\delta d$ restricted to $H=(\ker d)^\bot$ and the associated heat
$e^{-t\delta d}$ on $H$. Since these the second operator is the
Laplace transform of the former, the two asymptotics of
$\mathrm{Tr}_\Gamma (e^{-t\delta d})$ when $t\rightarrow +\infty$ and
of $F_{\delta d}(\lambda)= \tr_\Gamma (\Pi (\delta d \leq \lambda))$
when $\lambda \rightarrow 0$ are related in the following way (see
appendix of \cite{Gr-Sh}).  One has
$$
\tr_\Gamma (e^{-t\Delta_p}\Pi_H) \asymp t^{-\alpha_p/2}\ 
\mathrm{when}\ t\rightarrow +\infty
$$
iff
$$
F_{\delta d}(\lambda) \asymp \lambda^{\alpha_p/2}\ \mathrm{when}\ 
\lambda\rightarrow 0\,.
$$
A more general definition of $\alpha_p$, called the $p$th
Novikov-Shubin number of $M$, is
\begin{equation}
  \label{eq:1:1}
  \alpha_p = 2\liminf_{\lambda \rightarrow 0} \bigl(\ln F_{\delta
  d}(\lambda) / \ln \lambda\bigr) \in [0,+\infty]
\end{equation}
(also it seems there is no known geometric example where this liminf
is not an actual limit.)

The first one, $\alpha_0$, describing heat decay on functions is
known.
\begin{thm}\label{thm:1:2} (Varopoulos \cite{Va}, Gromov)
  $\alpha_0 < +\infty$ iff $\Gamma$ has polynomial growing that is iff
  $\Gamma$ is (virtually) nilpotent, and then
  $$\alpha_0= \mathrm{growing\ rate\ of\ } \Gamma = \lim_{R\rightarrow
    +\infty} \frac{\ln \mathrm{Vol_\Gamma} (B(0,R))}{\ln R}\in \N.$$
\end{thm}
It appears then that $\alpha_0$ is a rough invariant of $M$,
independent of the metric $g$, depending only here on the large scale
structure of $\Gamma=\pi_1(M)$. For the other $\alpha_p$, one has
\begin{thm}\label{thm:1:3}
  (Gromov-Shubin, Efremov \cite{Gr-Sh}) $\alpha_p$ are homotopy
  invariants of $M$.
\end{thm}
In fact it turns out that $\alpha_p$ depends on the (rational)
homotopy type of $M$ up to degree $p$. In particular the next one
$\alpha_1$, describing heat on $1$-forms, depends only on $\pi_1(M)$,
and is thus of particular interest. We will mainly focus on its study
here, although some sections will deal with related fields.

The (very partial) results are actually disseminated along the paper.
We will illustrate them with (carefully chosen) examples. This work
relies on the use of constructions that have already been exposed in
\cite{Ru}.  Anyway we recall them first.

\section{Some differential geometry on C-C manifolds}
\label{sec:2}

\subsection{Dilations and differential on graded groups}
\label{sec:2:1}

As time increases, heat spreads in larger and larger domains.  It is
therefore convenient to re-scale the phenomena. This requires some
dilation acting on the space. Typical groups admitting such a
structure are graded nilpotent (Lie) groups. These are groups $G$ such
that their Lie algebra $\mathfrak{g}$ splits in
$$
\mathfrak{g}=\bigoplus_{i=1}^r \mathfrak{g}_i \ \mathrm{with\ 
  }[\mathfrak{g}_1, \mathfrak{g}_i] \subset \mathfrak{g}_{i+1}.
$$
We call $r$ the \emph{rank} of $G$.  The multiplication by $i$ on $
\mathfrak{g}_i$ induces a family of dilations $h_\varepsilon$ on $G$
by $h_\varepsilon(\exp X_i) = \varepsilon^i (\exp X_i)$ for $X_i \in
\mathfrak{g}_i$. The function $w=i$ on ${\mathfrak g}_i$ is called the
\emph{weight} of vectors, and similarly on covectors. It can be
extended to the whole differential algebra of $G$. Forms of weight $p$
are spanned by
$$\theta^{i_1}\wedge \theta^{i_2}\wedge \cdots \wedge \theta^{i_k}
\mathrm{\ with \ }\sum_{j=1}^k w(\theta^{i_j}) =p.$$

De Rham's differential $d$ is not homogeneous with respect to the
weight. In fact it splits like
\begin{equation}
  \label{eq:2:1}
  d=d_0+d_1+\cdots +d_r
\end{equation}
with components $d_k$ increasing the weight by $k$. Indeed on
functions $d_0f=0$ since there is no $1$-forms of weight $0$, whereas
$d_kf$ for $1\leq k\leq r$ is just $df$ restricted to vectors of
weight $k$. Now, if $\alpha$ is a (left) invariant form of weight $p$
one has $d(f\alpha)=df\wedge \alpha + fd\alpha$. Then for $k\geq 1$,
$d_k(f\alpha) = d_kf\wedge \alpha$, while
$$
d_0(f\alpha) = fd_0 \alpha=fd\alpha
$$
is an \emph{algebraic} (order $0$) operator. The assertion that for
invariant forms $d\alpha=d_0\alpha$ is of the same weight as $\alpha$
is exactly dual to the property that $[\mathfrak{g}_i, \mathfrak{g}_j]
\subset \mathfrak{g}_{i+j}$, as is seen starting from invariant
$1$-forms for which $d\theta(X,Y)=-\theta ([X,Y])$ on invariant
vectors.

Now let $G$ be endowed with an invariant metric such that the
$\mathfrak{g}_i$ are mutually orthogonal.  If $\alpha$ is a form of
weight $p$, one has point-wise $||h_\varepsilon^*\alpha ||=
\varepsilon^p ||\alpha||$, so that
$$
||d(h_\varepsilon^*\alpha)||\geq \varepsilon^p||d_0 \alpha||= ||d_0
(h_\varepsilon^*\alpha)||.
$$
Therefore $d_0$, the zero order part of $d$, invisible on functions
and $\R^n$, is an obstruction of decreasing the differential of forms
through dilations. In other words, the spectral spaces $E(\delta d
\leq \lambda)$ have to contract on $\ker d_0$ in some sense when
$\lambda\rightarrow 0$.

\subsection{Cutting $d_0$ out of $d$}
\label{sec:2:2}

One would like to find a way of working directly on $\ker d_0$ when
dealing with small spectrum problems on our groups. The spectral
sequence technique achieves this from the algebraic viewpoint. In fact
the previous discussion pointed out a natural decreasing filtration of
$\Omega^*G$ by the spaces $F^p$ of forms of weight $\geq p$. The
decomposition \eqref{eq:2:1} says that $d$ respects it. Thus we have a
filtered complex giving rise to a sequence of spaces converging to the
(graded part of) de Rham cohomology of $G$.

We briefly recall this technique, $d_0$ is interpreted as $d$ acting
on $F^p/F^{p+1}$. Let $E_0= \ker d_0/ \im d_0$ denotes its cohomology.
This is the bundle which fiber is the Lie algebra cohomology
$H^*({\mathfrak g},\R)$. We get then $d$ acting on the quotient space
$E_0$, whose cohomology is called $E_1$, etc. Note that this process
produces successive quotients by image of \emph{differential}
operators, not very easy to handle analytically.  The fact that $d_0$
is algebraic will allow to perform part of the construction
\emph{inside} the de Rham complex. This will be useful to keep track
of the original problem of analyzing the small spectrum of $d$.

In fact the constructions apply to more general spaces than graded
groups called Carnot-Carath\'eodory manifolds which we introduce now.

\subsubsection{C-C manifolds and graded groups}
\label{sec:2:2:1}

\begin{defn}
  An equiregular Carnot-Carath\'eodory structure on a manifold $M$
  consists in an (increasing) filtration of $TM$ by bundles $H_i$ such
  that
  \begin{equation}
    \label{eq:2:2:1}
    [H_1,H_i] \subset H_{i+1} \ \mathrm{with}\ H_r =TM \ \mathrm{for\
    some\ }r.
  \end{equation}
\end{defn}
\begin{rems}
\item - Strictly speaking this definition is more general than the
  usual one given in \cite{Gr} for instance. Here we don't require
  $TM$ to be generated by $H_1$ through brackets, that is $H_1$ to be
  an H\"ormander distribution and $H_{i+1}=[H_1, H_i]$. This will
  allow us to change the grading if necessary.
\item - The equiregularity assumption means that we suppose the
  distributions $H_i$ are bundles. Their dimensions are not allowed to
  jump.
\end{rems}
The link with our previous discussion is the well known fact that an
equiregular C-C manifold admits a tangent graded nilpotent Lie algebra
${\mathfrak g}_{x_0}$ at each point. Indeed due to \eqref{eq:2:2:1}
the usual bracket on $TM$ admits a quotient map
$$
[\ , \ ]_0: H_p/H_{p-1} \times H_q/H_{q-1} \rightarrow
H_{p+q}/H_{p+q-1}.
$$
which is an algebraic operator (order $0$), since $[X,fY]_0 =
\Pi_{H_{p+q}/H_{p+q-1}}(f[X,Y] + (X.f))= f[X,Y]_0$. Therefore $[\ , \ 
]_0$ defines a graded Lie algebra structure on
\begin{equation}
  \label{eq:2:2:2}
  \mathfrak{g}_{x_0}= \bigoplus_{i=1}^r \mathfrak{g_i} \ 
\mathrm{with}\ \mathfrak{g_i} = H_i/H_{i-1}\ \mathrm{at\ }x_0.
\end{equation}
We recover also a decreasing filtration of $\Lambda^*T^*M$. Firstly
vectors in $H_p$ are called of weight $\leq p$, and dually a
differential form of degree $k$ is of weight $\geq p$ iff it vanishes
on all sets of $k$ vectors of total weight $< p$. Let $F_p$ be the
bundle of forms of weight $\geq p$. Again $F_p$ is $d$-stable as come
from \eqref{eq:2:2:1}, and we get a filtered complex on $\Omega^*M$.
As before $d_0=d$, acting on $F_p/F_{p-1}$, is an algebraic operator.
This zero order part of $d$ is identified with the previous $d_0$ we
introduced on ${\mathfrak g}_{x_0}$.

\begin{examples}\label{ex:2:3}
\item - Let $D^k$ be a distribution of $k$-planes in $\R^n$.
  Generically in the jets of the vector fields generating it, $D^k$
  will not be integrable and even will be an H\"ormander distribution.
  This gives a lot of (true) C-C structures on $\R^n$.
\item - We just describe some classical cases. An hyperplane
  distribution $D^{2n}$ in $\R^{2n+1}$ is called a contact structure
  when the Lie bracket $[\ ,\ ]_0$ is non-degenerate on $D$.  Reducing
  this bilinear form gives that $G_{x_0}=\exp{\mathfrak g}_{x_0} $ is
  isomorphic to the Heisenberg group. Even more, thanks to Darboux'
  theorem the contact structure itself is locally isomorphic to its
  tangent group $G_{x_0}$.
\item - A generic $D^2$ in $\R^4$, called an Engel's structure, also
  gives rise to a unique tangent $\mathfrak{g}_{x_0}$.  This is the
  $4$-dimensional Lie algebra of rank $3$ with the following relations
  $[X,Y]_0 = Z$, $[X,Z]_0 = T$ and $0$ elsewhere. Here
  $D=H_1=\mathrm{span}(X,Y)$, $H_2= [D,D]= \mathrm{span}(X,Y,Z)$, and
  $H_3 =\R^4 =[D,[D,D]]=\mathrm{span}(X,Y,Z,T)$. Again $D$ is an
  H\"ormander distribution, but we can look this structure
  differently. Consider $H'_1=\mathrm{span}(X)$, $H'_2=
  \mathrm{span}(X,Y)=D$, $H'_3 =H_2$, $H'_4 =H_3=\R^4$. This C-C
  structure has the same model group $G_{x_0}$ as before, but with a
  different grading. This will give another way of rescaling the
  asymptotic problems we consider.
\end{examples}

\subsubsection{One homotopy and three complexes}
\label{sec:2:2:2}

We now show how to retract the de Rham complex of a C-C manifold $M$
on its $E_0$ part. Recall that $E_0= \ker d_0/\im d_0=
H^*(\mathfrak{g}_{x_0})$, where $\mathfrak{g}_{x_0}$ is tangent Lie
algebra at $x_0$. In order to work on bundles we will assume now that
$\mathrm{dim} (E_0)$ is constant (always satisfied on an open dense
set of $M$).
\begin{defn}
  A C-C structure is $E_0$-regular if $E_0$ is a bundle.
\end{defn}

Choose a metric on $M$. Let $V_i$ be the orthogonal supplement of
$H_i$ in $H_{i-1}$. By \eqref{eq:2:2:2}, this fixes an isometry
between $T_{x_0}M=\bigoplus_{i=1}^r V_i$ and $\mathfrak{g}_{x_0}=
\bigoplus_{i=1}^r \mathfrak{g_i}$ and fixes the weight on $TM$. This
allows to see $d_0$ as acting between spaces of (true, not quotiented)
forms on $M$ of given pure weight. We identify $E_0$ with $\ker d_0
\cap \ker \delta_0=\ker d_0 \cap (\im d_0)^\bot$, where $\delta_0$ is
the metric adjoint of $d_0$. We also get a partial inverse of $d_0$ by
$d_0^{-1} = (\delta_0 d_0)^{-1} \delta_0$. We can use this $d_0^{-1}$
as a partial inverse of $d$.

Define the following retraction on $\Omega^*M$
\begin{equation}
  \label{eq:r}
  r=\id - d_0^{-1} d -d d_0^{-1}.
\end{equation}
This is an homotopical equivalence (preserving $d$ and the cohomology)
whose zero order term is $r_0=\id -d_0^{-1} d_ 0 -d_0 d_0^{-1} =
\Pi_{E_0}$ a projection, whereas $r-r_0= -d_0^{-1} (d-d_0)-
(d-d_0)d_0^{-1}$ strictly increases the weight and is therefore
nilpotent. This incites to iterate the homotopy $r$ in order to
retract the de Rham complex on the smallest possible space.  The maps
$r^k$ actually converge toward a map that certainly have to be both an
homotopical equivalence and a projection to a sub-complex along an
other. The following lemma is useful to identify quickly the limit
spaces and operators.
\begin{lemma}
  The map $d_0^{-1} d$ induces an isomorphism on $\im d_0^{-1}$, whose
  inverse is a differential operators $P$.
\end{lemma}
\begin{proof}(\cite{Ru})
  On $\im d_0^{-1}$, one can write $ d_0^{-1} d = \id + D$ where $D=
  d_0^{-1} (d-d_0)$ is nilpotent since it strictly increases the
  weight of forms. One has then
  $$
  P= (d_0^{-1} d)^{-1} = \sum_{k=0}^{d(M)}(-1)^k D^k
  $$
  where $d(M)= \mathrm{weight}(\Lambda^{\mathrm{max}}T^*M)$.
\end{proof}
Define $Q=Pd_0^{-1}$ with $P$ as in the lemma. One has then.
\begin{thm}
  \label{thm:2:2:1}
  Let $M$ be a $E_0$-regular C-C space.
  \begin{enumerate}
  \item $\Omega^*M$ splits in the direct sum of two sub-complexes
      \begin{displaymath}
        E= \ker d_0^{-1} \cap \ker (d_0^{-1} d) \ \mathrm{and}\ F=\im
        d_0^{-1} + \im (dd_0^{-1}). 
      \end{displaymath}
      The projection $\Pi_F$ on $F$ along $E$ is given by the
      differential operator $Qd+dQ$, with $Q$ as above.
    \item The homotopies $r^k$ converge to $\Pi_E$ the projection on
      $E$ along $F$. It is the homotopical equivalence given by the
      differential operator $\id -\Pi_F = \id -Qd -dQ$.
    \item One has
    \begin{equation}
      \label{eq:2:2:3}
      \Pi_{E_0} \Pi_E \Pi_{E_0} = \Pi_{E_0} \ \mathrm{and}\ \Pi_E
      \Pi_{E_0}\Pi_E = \Pi_E,
    \end{equation}
    saying that $E$ and $E_0$ are in bijection, and that $\Pi_E$
    restricted to $E_0$ and $\Pi_{E_0}$ restricted to $E$ are inverse
    maps of each other. In particular the complex $(E,d)$ is
    conjugated to another one $(E_0,d_c)$ with $d_c= \Pi_{E_0} d\Pi_E
    \Pi_{E_0}$.
  \end{enumerate}
\end{thm}
The constructions are summarized in the diagram
\begin{displaymath}
  \xymatrix{
    F \ar@{^{(}->}@<1ex>[r]^i \ar[d]_d & \Omega^*M \ar@<1ex>[l]^{\Pi_F}
    \ar@<1ex>[r]^{\Pi_E} \ar[d]_d & E \ar@{_{(}->}@<1ex>[l]^i 
    \ar@<1ex>[r]^{\Pi_{E_0}} \ar[d]_d & E_0
    \ar@<1ex>[l]^{\Pi_E}\ar[d]^{d_c}  \\
    F \ar@{^{(}->}@<1ex>[r]^i & \Omega^*M \ar@<1ex>[l]^{\Pi_F}
    \ar@<1ex>[r]^{\Pi_E} & E \ar@{_{(}->}@<1ex>[l]^i
    \ar@<1ex>[r]^{\Pi_{E_0}} & E_0 \ar@<1ex>[l]^{\Pi_E}
    }
\end{displaymath}
A short proof is given in \cite{Ru}(thm $1$). We don't repeat it here,
but insist on particular points.

Observe that $E$ is a space of forms satisfying some differential
equations ($E=\ker d_0^{-1} \cap \ker ( d_0^{-1} d)$) that projects
bijectively onto the algebraic $E_0 = \ker d_0 \cap \ker d_0^{-1}$.
Even more, the equation $\Pi_{E_0} \Pi_E \Pi_{E_0} = \Pi_{E_0}$ says
that $\Pi_E$ restricted to $E_0$ is $\id + $a $(E_0)^\bot$ part. In
other words $E$ is a particular space of liftings (extensions) of
$E_0$. For computations, we notice from \eqref{eq:r} that the
retraction $r$ preserves the space $(\ker d_0)^{-1} = \ker \delta_0
\supset E_0 $ on which
\begin{displaymath}
  r = \id - d_0^{-1}(d-d_0) 
\end{displaymath}
with $d_0^{-1}(d-d_0)$ strictly increasing the weight. Starting from
some form $\alpha$ of weight $p$ in $E_0$ this gives by iteration the
successive extensions in greater weight by
\begin{equation}
  \label{eq:2:2:5}
  \left\{
  \begin{split}
    (\Pi_E \alpha)_p &= \alpha \\
    (\Pi_E \alpha)_{p+k+1} &= -d_0^{-1} \Bigl(\sum_{l=1}^r d_l (\Pi_E
    \alpha)_{p+k+1-l}\Bigr)
  \end{split}\right.
\end{equation}
where $d_l$ is the part of $d$ that increases the weight by $l$.
\begin{rems}
\item - At this point it is clear that this construction is everything
  but a new idea, merely an hold one in a fancy dress. Namely this
  $\Pi_E$ is a realization of the first homotopical equivalence that
  is predicted by the general fact that the spectral sequence starting
  on $E_0$ will finally compute (the graded part of) the cohomology.
  The map $d_c$ describes the remaining obstructions of extending a
  form in $E_0$ to a true closed one, put all together.
\item - Formally at least it is tentative to carry on with the
  retraction process by replacing $d_0$ by $d_{c,1}$, the part of
  $d_c$ that increases the weight by one. Choosing a supplement space
  to $\ker d_{c,1}$ as the $L^2$ closure of $\im d_{c,1}$, leads to a
  (now not bounded) partial inverse $d_{c,1}^{-1}$ of $d_{c,1}$ and an
  another retraction $r_1= \id - d_{c,1}^{-1} d_c - d_c d_{c,1}^{-1}$.
  It should be iterated and so on.  Of course this does not seem
  realistic analytically. Anyway some analysis (yet a little bit
  mysterious) related to the spectral sequence structure is working.
  We will describe it in section \ref{sec:5}.
\item - We didn't take care of the invariance of $d_c$ on the choice
  of metric, mainly because it is not invariant in general ! This
  doesn't matter here.  Again this construction is taken as a
  convenient approximation of the underlying invariant spectral
  sequence, which is hopefully related to the (even homotopically)
  invariant problem of studying the asymptotic heat decays.
\end{rems}

We conclude this general section with a remark on duality. As observed
in \cite{Ru}, the complex $(E_0, d_c)$ is Hodge $*$-dual. This should
look rather surprising since the retraction map $r=\id - d_0^{-1} d -d
d_0^{-1}$ breaks the symmetry between $d_0$ and $\delta_0$. Anyway, we
have
\begin{prop}
  \begin{enumerate}
  \item $* \delta_0 = (-1)^{k+1} d_0*$ on $E_0^k$ and $*$ preserves
    $E_0$.
  \item $*E$ is orthogonal to $F$, equivalently the pairing $(\alpha,
    \beta) \rightarrow \int_M \alpha \wedge \beta$ vanishes on
    $E\times F$. One has $*\Pi_E'=\Pi_E *$, where $\Pi_E'$ is the
    formal adjoint of $\Pi_E$.
  \item $*\delta_E = (-1)^{k+1} d_E *$ on $k$-forms. Similarly
    $*\delta_c = (-1)^{k+1} d_c*$ on $E_0^k$.
  \end{enumerate}
\end{prop}
\begin{proof}
  - The first observation comes from the same (point-wise) duality
  standing at the Lie algebra level on $\mathfrak{g}_{x_0}$.
  
  - By definition $F= \im \delta_0 + \im d\delta_0$ and therefore
  $$
  *(F^\bot) =*( \ker d_0 \cap \ker d_0 \delta) = \ker \delta_0 \cap
  \ker d\delta_0 =E.
  $$
  Then,
  \begin{align*}
    (\alpha, \Pi_E' \beta) &= (\Pi_E \alpha , \beta) = \int_M \Pi_E
    \alpha \wedge *\beta \\
    & = \int_M\Pi_E \alpha \wedge \Pi_E *\beta = \int_M \alpha \wedge
    \Pi_E *\beta \\
    & = (\alpha, *^{-1} \Pi_E *\beta).
  \end{align*}
  
  - From $d_E = d\Pi_E =\Pi_E d$, we find that
  \begin{displaymath}
    *\delta_E = *\Pi_E' \delta = \Pi_E *\delta= (-1)^{k+1} \Pi_E d * =
     (-1)^{k+1} d_E *,
  \end{displaymath}
  and similarly for $\delta_c$ starting from $d_c= \Pi_{E_0} d_E
  \Pi_{E_0}$ and using that $\Pi_{E_0}$ commutes with $*$ and is
  self-adjoint, being an orthogonal projection.
\end{proof}

\subsection{A few examples}
\label{sec:2:3}

We briefly describe some specific cases of the previous constructions.

\item $\bullet$ On $\R^n$, or $M^n$ with the trivial C-C structure
  ($H_1 = TM$), one has $d_0=0$, and $d_c=d$.
  
\item $\bullet$ Let $(M^{2n+1},H)$ be a contact structure. Let
  $\theta$ be such that $H=\ker \theta$. $\Omega^*M$ splits in
  $\Omega^*H + \theta \wedge \Omega^*H$ and the $0$ order part of
  $d(\alpha + \theta \wedge \beta)$ is readily seen to be $d\theta
  \wedge \beta$. From this we get that $E_0^k$ are primitive
  horizontal $k$-forms if $k\leq n$, and the vertical co-primitive
  ones ($\in \ker (d\theta \wedge\id )$) if $k\geq n+1$. These bundles
  $E_0^k$ are of pure C-C weight, namely $k$ if $k \leq n$ and $k+1$
  if $k\geq n+1$.
  
  This implies that the $d_c$ complex consists in first order
  operators, except an order $2$ one in degree $n$, due to the
  ``jump'' of weights between $E_0^n$ (weight $n$) and $E_0^{n+1}$
  (weight $n+2$). This is the contact complex.
  
\item $\bullet$ Let $M=G$ be the Engel group (see examples
  \ref{ex:2:3}). Recall this is the $4$ dimensional Lie group, with
  relations $[X,Y]=Z$, $[X,Z] =T$, and other brackets vanishing. In
  the dual base of $1$-forms, this translates in
  $$d_0 \theta_X = d_0 \theta_Y = 0,\quad d_0 \theta_Z = -\theta_X
  \wedge \theta_Y \quad \mathrm{and}\quad d_0 \theta_T = -\theta_X
  \wedge \theta_Z.$$
  We find then
  \begin{equation*}
    \left\{
      \begin{split}
        E_0^0 &=C(G)\ :\ \mathrm{functions \ on\  } G,\\
        E_0^1 &= \mathrm{span} \left(\theta_X,\theta_Y\right)\ :\ 
        \mathrm{horizontal\ one\ forms}, \\
        E_0^2 &= \mathrm{span} \left(\theta_Y \wedge \theta_Z,
          \theta_X
          \wedge \theta_T\right) = * E_0^2,\\
        E_0^3 &= *E_0^1= \mathrm{span}\left( \theta_Y \wedge \theta_Z
          \wedge \theta_T, \theta_X \wedge \theta_Z \wedge
          \theta_T\right) \\
        E_0^4 &= \Omega^4G.
      \end{split}
    \right.
  \end{equation*}
  We compute $d_c$. As there is nothing to lift on functions $\Pi_E
  f=f$, and
  \begin{displaymath}
    d_cf = \Pi_{E_0^1}df = d_{H_1}f = (X.f)\theta_X + (Y.f)
    \theta_Y.
  \end{displaymath}
  Using \eqref{eq:2:2:5}, we can compute $\Pi_E \alpha =
  \tilde{\alpha}$ for $\alpha = \alpha(X) \theta_X + \alpha (Y)
  \theta_Y \in E_0^1$.  One gets
  \begin{align*}
    - d_0 (\tilde{\alpha} (Z) \theta_Z) & = \tilde{\alpha}
    (Z) \theta_X \wedge \theta_Y \\
    &= d_1 \alpha = (X.\alpha (Y) -Y.\alpha (X)) \theta_X \wedge
    \theta_Y,
  \end{align*}
  thus $\tilde{\alpha} (Z) = X.\alpha (Y) -Y.\alpha (X)$. Again with
  \eqref{eq:2:2:5},
  \begin{align*}
    - d_0 (\tilde{\alpha} (T) \theta_T) &= \tilde{\alpha}(T)
    \theta_X \wedge \theta_Z \\
    & = (X.\tilde{\alpha} (Z) - Z. \alpha(X)) \theta_X \wedge
    \theta_Z,
  \end{align*}
  and $\tilde{\alpha} (T) = X.\tilde{\alpha} (Z) -Z.\alpha (X) =
  X^2.\alpha (Y) - (XY+Z).\alpha (X)$. Now $d_c \alpha$ is the
  restriction of $d\tilde{\alpha}$ to $E_0^2$, that is
  \begin{displaymath}
    d_c \alpha = (Y.\tilde{\alpha} (Z) -Z. \alpha(Y)) \theta_Y \wedge
    \theta_Z + (X.\tilde{\alpha}(T) -T. \theta(X)) \theta_X \wedge
    \theta_T. 
  \end{displaymath}
  The full complex $d_c$ may be completed either by $*$-duality or
  computing $\Pi_E$ in degree $2$.  The result can be read out from
  the following diagram, adding all possible travels between points
  gives the various components of the liftings and $d_c$
  \begin{equation}
    \label{eq:2:2:6}
    \xymatrix@=6pt{
      & (\theta_Y) \ar[rr]^{-Z} \ar[dr]^X &&(\theta_{Y \wedge Z})
      \ar[dr]_{-X} \ar[drrr]^T && && \\
      (f) \ar[ur]^Y \ar[dr]_X & & \theta_Z \ar[ur]_Y \ar[dr]^X
      && \theta_{Y \wedge T} \ar[dr]_X \ar[rr]_{-Z} && (\theta_{Y\wedge Z
      \wedge T}) \ar[dr]^X &\\
      &(\theta_X) \ar[ur]^{-Y} \ar[rr]^{-Z} \ar[drrr]_{-T} && \theta_T
      \ar[dr]^X && \theta_{Z \wedge T} \ar[ur]_Y \ar[dr]^X && (\Omega^4G)\\
      &&&& (\theta_{X \wedge T}) \ar[ur]^{-Y} \ar[rr]^{-Z} &&
      (\theta_{X\wedge Z\wedge T}) \ar[ur]^{-Y} &
      }
  \end{equation}
  
  It should be noted that $d_c$ doesn't depend on the various possible
  choices of grading on $G$ (see \ref{ex:2:3}), if keeping the same
  basis for $E_0= \ker d_0 / \im d_0$.
  
  Now, what about $d_c$ on a general Engel structure $D^2 \subset
  \R^4$ (or $TM^4$)? Unlike the contact case these are not locally
  diffeomorphic to their tangent group $G$. That means that we can't
  find a local system of vector fields $X,Y$ generating $D^2$ and
  satisfying \emph{exactly} the previous bracket relations. Even so,
  thanks to $[\ ,\ ]_0$, they can be satisfied up to vectors of lower
  C-C weight. The conclusion is that the general $d_c$ will be a
  perturbation of the above $d_c$ on $G$ by differential operators of
  lower C-C weight. This is of course a general feature of the
  construction.

\item $\bullet$ We close this series of examples with a glimpse toward
  the nice case of C-C structures given by a generic $4$ dimensional
  distribution $D^4$ in $\R^7$ (or $TM^7$). Here again we are in an
  exceptional situation where there are finitely many isomorphic type
  of possible tangent graded group $G_{x_0}$.
  
  We see this. Its Lie algebra $\mathfrak{g}_{x_0}$ is generated by
  $\mathfrak{g}_1 = D$ and $\mathfrak{g}_2 = TM/D$ and determined by
  its curvature $d_0 : \Lambda^1 (TM/D)^* \rightarrow \Lambda^2
  D^{4*}$ given by $d_0 \theta (X,Y) = -\theta ([X,Y]_0)$. This map is
  injective iff $[\ ,\ ]_0: \Lambda^2D \rightarrow TM/D$ is
  surjective, that is the distribution $D$ is bracket generating. In
  that (generic) case $L= \im d_0$ is a $3$ dimensional subspace in
  $\Lambda^2D^*$. This is a famous case where $L$ is determined up to
  isomorphism by the signature of the quadratic form $q(\omega)
  d\mathrm{vol}_{D} = \omega \wedge \omega$ restricted to it. The
  distribution $D$ is called elliptic if $q$ is positive definite on
  $L$ (changing the orientation of $D$ if necessary). This is an open
  condition, but not dense due to the other open possibility of an
  hyperbolic $(2,1)$ signature !
  
  If $D$ is elliptic, our assumption by now, one knows there exists a
  unique conformal class of metric $g$ on $D$ such that $L$ becomes
  $\Lambda^{2,+} D^*$ the space of $*$-self-dual $2$ forms of $D^4$.
  It is convenient to describe $L$ by seeing $D$ as the quaternions
  $\H$ in which case $L= \mathrm{span}\left\{ d\theta_l\right\}_{1\leq
    l\leq 3}$ is such that
  $$
  d_0\theta_l = g(J_l\cdot, \cdot),
  $$
  for $J_l=i,j,k \in \H$. This is an example of an Heisenberg type
  group or $H$-group (see \cite{Ka}, \cite{Co-Do}), associated here to
  the quaternion-hyperbolic group $Sp(2,1)$ .
  
  The cohomology $E_0$ of this structure is easily computed.
  \begin{enumerate}
  \item Again $E_0^1 = \Lambda^1 D^*$ are the horizontal $1$-forms.
  \item $E_0^2$ splits in two weights $E_0^{2,(2)} = \Lambda^2 D^*/ L
    \simeq \Lambda^{2,-}D^*$, the anti self-dual part of
    $\Lambda^2D^*$. The weight $3$ part $E_0^{2,(3)}$ is generated by
    $\theta_1 \wedge J_1 \alpha -\theta_2 \wedge J_2 \alpha$ and
    $\theta_2 \wedge J_2 \alpha -\theta_3 \wedge J_3 \alpha$ for
    $\alpha \in \Lambda^1D^*$.
  \item $E_0^3$ is seen to be of pure weight $4$. It is the
    $14$-dimensional space of $3$-forms $\gamma = \sum_{i=1}^3\theta_i
    \wedge \beta_i$ with $\beta_i \in \Lambda^2D^*$ satisfying $d_0
    \gamma = \sum_{1\leq i\leq 3}\left(d\theta_i ,\beta_i
    \right)\mathrm{dvol}_D =0$ modulo $\im d_0$ generated by $d_0
    (\theta_1 \wedge \theta_2) = - \theta_1 \wedge d\theta_2 +
    \theta_2 \wedge d\theta_1$ and the two others coming from
    permutation.
  \item The missing degrees are obtained by duality.
  \end{enumerate}
  
  We summarize the information in the following diagram, with the
  degree of forms along x and $(\mathrm{weight} - \mathrm{degree})$
  along y
  \begin{equation}
    \label{eq:2:2:7}
    \xymatrix@=7pt{
      &&&&& E_0^{5,(8)} \ar[r]^{d_c} & E_0^{6,(9)} \ar[r]^{d_c} &
      \Omega^7M \\
      &&&& E_0^{4,(6)} \ar[r]^{d_c} \ar[ur]^{d_c} & E_0^{5,(7)}
      \ar[ur]_{d_c} &&\\
      && E_0^{2,(3)} \ar[r]^{d_c} &  E_0^{3,(4)} \ar[ur]^{d_c}&&&&  \\
      f \ar[r]^{d_D} & \Lambda^1D^* \ar[r]^{d^{-}} \ar[ur]^{d_c} &
      \Lambda^{2,-}D^* \ar[ur]_{d_c} &&&&&
      }
  \end{equation}
  The interesting feature here is that the bottom line looks like an
  elliptic complex on the $4$ dimensional $D$, except that $D$ is not
  integrable here. Anyway the ellipticity of this first order part of
  $d_c$ will be helpful in in determining the heat decays of forms on
  the tangent group.

\subsection{Geometric comments in degree $\leq 2$}
\label{sec:2:4}

We give some precisions on the structure of $d_c$ in degree $\leq 2$.

\item $\bullet$ We start with $E_0^1=H^1(\mathfrak{g}_{x_0})$ on C-C
  manifold $M$. If $D$ is bracket generating, that is an H\"ormander
  distribution in $TM$, then $\mathfrak{g}_{x_0}$ will be generated by
  its vectors of weight $1$. (These particular graded groups are often
  called filtered groups.) This is equivalent to the fact that $d_0$
  is injective on $1$ forms of weight $\geq 1$. Indeed $d_0$
  restricted to one forms is dual to the bracket map $[\ , \ ]_0 :
  \mathfrak{g}_{x_0} \wedge \mathfrak{g}_{x_0} \rightarrow
  \mathfrak{g}_{x_0}$ and this one is surjective onto vectors of
  weight $\geq 1$ precisely when $D$ is bracket generating.
  
  As a consequence, $E_0^1=\ker d_0 /\im d_0$ consists only of forms
  of weight $1$, which can be identified with $\Lambda^1D^*$, the
  partial $1$-forms on $D$. Then $d_c f$ is simply the restriction of
  $d f$ to vectors in $D$, a first order operator.
  
\item $\bullet$ Assume again that $\mathfrak{g}_{x_0}$ is filtered,
  that is $D=\mathfrak{g}_1$ is bracket generating. We now see how
  $E_0^2 = H^2(\mathfrak{g}_{x_0})$ is linked to the relations
  defining $G_{x_0}$.  Let $\widetilde{G}$ be the free Lie group
  generated over $D$. The map $\Pi: \widetilde{G} \rightarrow G_{x_0}$
  is surjective since $D$ is bracket generating.  Therefore one has
  $G_{x_0} = \widetilde{G}/N$ where $N =\Pi^{-1}(0)$ interprets as the
  normal subgroup of relations of $G$ (with respect to
  $\widetilde{G}$). At the Lie algebra level $N=\exp \mathfrak{n}$
  with $\mathfrak{n} = \ker \pi$ is the ideal of relations of
  $\mathfrak{g}_{x_0}$.
  
  Since $\mathfrak{n}$ is an ideal, $r = \mathfrak{n}/ [\mathfrak{n},
  \widetilde{\mathfrak{g}}]$ may be viewed as its space of generators,
  isomorphic to $R=N/(N,\widetilde{G})$ generating the relations $N$
  of $G_{x_0}$. Now a (classical) fact is that
  \begin{prop}\label{prop:2:9}
    $\Lambda^1 r^*$ is naturally isomorphic to
    $H^2(\mathfrak{g}_{x_0})$.
  \end{prop}
  \begin{proof}
    This is an homological consequence of the short exact sequence
    $\mathfrak{n} \rightarrow \widetilde{\mathfrak{g}} \rightarrow
    \mathfrak{g}_{x_0}$. We give the principle in our case. Consider
    $d_0:\Lambda^1 \widetilde{\mathfrak{g}}^* \rightarrow \Lambda^2
    \widetilde{\mathfrak{g}}^*$ given by $d_0 \alpha (X,Y) =-\alpha
    ([X,Y])$.
  
    The space $\Lambda^1 r^*$ identifies with $1$-forms on
    $\mathfrak{n}$ vanishing on $[\mathfrak{n},
    \widetilde{\mathfrak{g}}]$.  Let $\alpha \in \Lambda^1 r^*$, and
    $\overline{ \alpha}$ be any extension of $\alpha$ to $\Lambda^1
    \widetilde{\mathfrak{g}}^*$. Then $d_0 \overline{\alpha}$ vanishes
    on $\mathfrak{n} \times \widetilde{\mathfrak{g}}$ and is therefore
    the pull-back of a $2$-form $\beta$ on $\mathfrak{g}_{x_0}$. One
    checks easily that $\beta \in \ker d_{0,\mathfrak{g}_{x_0}}$, and
    that $\beta \rightarrow \beta + d_0 \gamma$ when the extension
    $\overline{\alpha}$ is changed (by a one form vanishing on
    $\mathfrak{n}$). We get therefore a map
    \begin{align*}
      [d_0] : \Lambda^1r^* & \rightarrow H^2 (\mathfrak{g}_{x_0}) =
      \ker d_0 /\im
      d_0 \\
      \alpha & \rightarrow [\beta = d_{0,\widetilde{\mathfrak{g}}}
      \overline{\alpha}].
    \end{align*}
     
    - Injectivity of $[d_0]$. If $ [d_0] \alpha =0$ then $\exists
    \gamma \in \Lambda^1\mathfrak{g}_{x_0}$ such that $d_0
    (\overline{\alpha} - \Pi^* \gamma) =0$. But since $D$ generates
    $\widetilde{\mathfrak{g}}$, we knows that $\overline{\alpha} -
    \Pi^* \gamma$ is a form of weight $1$. Therefore by restriction on
    ${\mathfrak n}$, of weight $\geq 2$, we get $\alpha =0$.
    
    - The surjectivity of $[d_0]$ relies on the fact that $H^2
    (\widetilde{\mathfrak{g}}) =0$.
    
    This can been taken as a definition of $\widetilde{\mathfrak{g}}$,
    meaning that there is one construction of
    $\widetilde{\mathfrak{g}}$ which is precisely done to achieve this
    ! (Starting with $D$, vectors of weight $1$, the space of
    covectors of weight $k+1$ in $\widetilde{\mathfrak{g}}$ is taken
    to be $\ker d_0$ among the \emph{two} forms of weight $\leq k+1$
    already constructed.) More geometrically, to $\beta \in \Lambda^2
    \widetilde {\mathfrak{g}}^* \cap \ker d_0$, a two cocycle,
    corresponds an extension of $\widetilde{\mathfrak{g}}$ by $\R$,
    $\mathfrak{g}'= \widetilde{\mathfrak{g}} \oplus \R \beta$ with $[\ 
    ,\ ]'$ defined by $[(x_1,t_1), (x_2, t_2)]' = \left( [x_1,x_2],
      \beta(x_1, x_2)\right)$. By universality of
    $\widetilde{\mathfrak{g}}$, this extension has to be trivial, a
    product, which corresponds to a $d_0$-exact $\beta$.
    
    We finish the proof. If $\beta$ is a $d_0$-closed two form of
    $\mathfrak{g}_{x_0}$, its pull-back $\Pi^*\beta$ on
    $\widetilde{\mathfrak{g}}$ is still $d_0$-closed and therefore
    $d_0$-exact.  Let $\alpha \in \Lambda^1
    \widetilde{\mathfrak{g}}^*$ be such that $\Pi^*\beta = d_0
    \alpha$.  One checks easily that the restriction of $\alpha$ on
    $\mathfrak {n}$ depends only on the cohomology class of $\beta$,
    and vanishes on $[\mathfrak{n},\widetilde{\mathfrak{g}}]$. It is
    therefore a one form on $r= \mathfrak {n}/ [\mathfrak{n},
    \widetilde{\mathfrak{g}} ]$.
\end{proof}

The previous proposition will give a geometric interpretation, in
terms of relations of a group $G$, of the pinching of heat decay on
$1$-forms that will be first related to the algebraic $E_0^2 =
H^2(\mathfrak{g})$. The use of the free Lie algebra $\widetilde
{\mathfrak{g}}$ is also helpful in understanding $d_c$ itself on
$1$-forms on a fixed filtered group $G$.
  
We see this. Let $\Pi: \widetilde{G} \rightarrow G$ as before and
$\alpha \in E_0^1(G)$. One has $\Pi^*\alpha \in E_0^1(\widetilde{G})$
and since $E_0^2(\widetilde{G})$ vanishes there will be no obstruction
in extending $\alpha$ on vectors of increasing weights to a closed one
$\beta$.  By injectivity of $d_0$ on $\widetilde{\mathfrak{g}}_{\geq
  2}$ this extension is unique. So one must have $\beta = \Pi_E
\Pi^*\alpha$ as given by the iteration of the retraction map $r$ on
$\widetilde{G}$ by \eqref{eq:r}. More concretely $\beta$ is determined
by
\begin{equation}
  \label{eq:2:4:1}
  \beta ([X_1, X_2]) = X_1. \beta (X_2) - X_2 . \beta (X_1)\quad
  \mathrm{and}\quad \beta = \Pi^*\alpha\ \mathrm{on}\ \mathfrak{g}_1.
\end{equation}
We note that in this extension process the components of $\beta$ stay
invariant functions along $\mathfrak{n}$
\begin{equation}\label{eq:2:4:2}
  Y.\beta(X) =0
\end{equation}
for $Y \in \mathfrak{n}$, as proved from \eqref{eq:2:4:1} by
recurrence on weight of $X$ and using that $\mathfrak{n}$ is an ideal.
Observe that the value of $\beta$ on $\mathfrak{n}$ is determined by
its value on the generators of $\mathfrak{n}$ because for any $Y \in
\mathfrak{n}$
\begin{equation}
  \label{eq:2:4:3}
  \beta ([X,Y]) = X. \beta (Y).
\end{equation}
This $\beta$ is a canonical form, but in order to push its restriction
to $\mathfrak{n}$ on $E_0^2(G)$ with proposition \ref{prop:2:9}, we
need to choose a supplementary subspace of $[\mathfrak{n},
\widetilde{\mathfrak g}]$ in $\mathfrak{n}$, or equivalently a set of
generators of relations of $G$. Depending on this choice, we finally
get $d_c \alpha$.

\section{Near-cohomology}
\label{sec:3}

So far we have only some formal hint in section \ref{sec:2:1} that the
asymptotic spectral problem we consider has something to do with the
previous constructions. The notion of near-cohomology, introduced in
\cite{Gr-Sh}, is the main tool that allows to relate explicitly the
problems.

\subsection{Definition and first applications}
\label{sec:3:1}

Here is a brief presentation
%``restricted overview''
of this notion (see \cite{Gr-Sh, Gr-Sh'} for details).  Let $M$ be a
complete riemannian manifold. An Hilbert complex over $M$, consists in
a sequence of Hilbert spaces $E$ over $M$ together with $d:D(d)
\subset E \rightarrow E$ such that $d$ are closed densely defined
operators and $d^2=0$. On the quotiented Hilbert space $F=E/\ker d$,
one considers the family of closed cones
\begin{displaymath}
  C_\varepsilon = \{ \alpha \in F, \|d\alpha\|_E \leq \varepsilon
  \|\alpha\|_F \}
\end{displaymath}
depending of $\varepsilon > 0$. These are shrinking (towards $\{0\}$)
as $\varepsilon \rightarrow 0$. By definition, the near-cohomology of
$(E,d)$ consists in this family of cones up to equivalence induced by
dilatational changes of $\varepsilon \rightarrow K \varepsilon$. That
means that two Hilbert complexes $(E,d)$ and $(E',d')$ have the same
near cohomology if for $\varepsilon$ small enough there exists a
constant $K> 0$ and \emph{bounded} maps $f:C_\varepsilon \rightarrow
C_{K\varepsilon}'$ and $g: C'_\varepsilon \rightarrow
C_{K\varepsilon}$ invertible on their images.

One observation of Gromov and Shubin is that bounded homotopical
equivalence between two Hilbert complexes induces equivalence of their
near-cohomologies.
\begin{defn}
  Two Hilbert complexes $(E,d)$, $(E',d')$ are homotopy equivalent, if
  there exists bounded maps $f: E \rightarrow E'$ and $g:E'\rightarrow
  E$ such that $fd = d'f$ on $D(d)$, $gd'=dg$ on $D(d')$, $g\circ f=
  \id_E +dA +B d$ on $D(d)$, resp. $f\circ g =\id_F +d'A' +B'd'$ on
  $D(d')$ for bounded operators $A,A',B,B'$.
\end{defn}
\begin{thm}\label{thm:3:2}
  (prop 4.1 \cite{Gr-Sh}) Let $(E,d)$, $(E',d')$ be two homotopy
  equivalent Hilbert complexes. Then, for $\varepsilon$ small enough,
  $f$ and $g$ induce maps $[f]: C_\varepsilon \rightarrow
  C_{K\varepsilon}' $ and $[g] : C'_\varepsilon \rightarrow
  C_{K\varepsilon} $ invertible on their images.
\end{thm}
\begin{proof}
  We recall the proof as it is short. Let $\alpha \in C_\varepsilon$.
  Then certainly
  \begin{equation}\label{eq:9}
    \|d'(f\alpha)\| = \|fd\alpha\| \leq \|f\| \|d\alpha\| \leq
    \varepsilon\|f\| \|\alpha\|.
  \end{equation}
  We need to control $\alpha$ by $\overline{f\alpha}$, the projection
  of $f\alpha$ in $F=E/\ker d'\simeq (\ker d')^\bot$. One has
  $\overline{f\alpha} = f\alpha + \beta$ with $\beta \in \ker d'$, so
  that
  \begin{displaymath}
    g(\overline{f\alpha}) = g(f\alpha) + g\beta = \alpha +d A\alpha
    +Bd\alpha  + g\beta,
  \end{displaymath}
  with $d A\alpha + g\beta \in \ker d$. Therefore, as $\alpha \in
  (\ker d)^\bot$,
  \begin{displaymath}
    \|\alpha\| \leq \|\alpha + d A\alpha + g\beta \| \leq
    \|g\overline{f\alpha} \| + \|B\| \varepsilon \|\alpha\|. 
  \end{displaymath}
  Finally
  \begin{displaymath}
    \|\alpha\| \leq \frac{\|g\|}{1 -\varepsilon \|B\|}
    \|\overline{f\alpha}\|, 
  \end{displaymath}
  giving the injectivity of $\overline{f}$ acting on $C_\varepsilon$,
  and the result together with \eqref{eq:9}.
\end{proof}

One can't apply directly this theorem to the complexes $(E,d)$ and
$(E_0,d_c)$ described in theorem \ref{thm:2:2:1}. This is because one
of the relevant map here, $\Pi_E$, being a differential operator, is
not bounded in $L^2$. Yet we observe that its un-boundedness occurs on
high energy forms, and we are precisely interested in the bottom of
the spectrum.  So we are leaded to consider an intermediate cut-off de
Rham complex $(E(\Delta \leq 1), d)$, where $E(\Delta \leq 1)$ is the
spectral space associated to $[0,1]$ by the Laplacian $\Delta$.

To use this remark, we first notice that the previous theorem
\ref{thm:3:2} applies to the full de Rham complex and the cut-off
one's, because $\Pi_{E(\Delta \leq 1)}= \id -dA -Ad$ where $A=\delta
\Delta^{-1}\Pi_{E(\Delta > 1)}$ is bounded. Thus, they have equivalent
near-cohomology (observe that although they have exactly the same
small spectrum their $C_\varepsilon$ are distinct). But now by
ellipticity of $\Delta$, the spectral projection $\Pi_{E(\Delta \leq
  1)}$ is a smoothing operator and therefore the differential operator
$\Pi_E =\id -Qd-dQ$ and $\Pi_{E_0}\Pi_E$ becomes bounded on $E(\Delta
\leq 1)$, making theorem \ref{thm:3:2} usable here.
\begin{thm}\label{thm:3:3}
  \cite{Ru} Let $M$ be a complete $E_0$-regular C-C manifold. Then the
  de Rham complex, $(E,d)$ and $(E_0,d_c)$ have equivalent
  near-cohomologies.
\end{thm}
\begin{proof}
  To complete the previous discussion, we make short comments on the
  closures of $(E_0, d_c)$ and $(E,d)$.
  
  - Firstly, $d_c$, being a differential operator, its formal adjoint
  $\delta_c$ has dense (initial) domain $C^\infty_0$ in $L^2$.
  Therefore $d_c$ is closable, for instance by $\overline{d_c} =
  (d_c)^{**}$ (see eg \cite{Re-Si}).
  
  - About $(E,d)$. Given a closed extension $\overline{d}$ of de
  Rham's differential, its restriction to ${\overline E}= \ker
  d_0^{-1} \cap \ker d_0^{-1}{\overline d}$ may be seen as a closed
  extension of $(E,d)$.
\end{proof}

The boundedness of differential operators on $E(\Delta \leq 1)$ has
also a geometric consequence on the shape of the cones $C_\varepsilon
(E(\Delta \leq 1),d)$.

We denote by $[\Pi_E]$ and $[\Pi_F]$ the actions of the retractions
$\Pi_E$ and $\Pi_F$ of theorem \ref{thm:2:2:1} on the quotient space
$\Omega^*M /\ker d$. One has still $[\Pi_F] =\id - [\Pi_E]$ with
$[\Pi_E]$ a projection, so that they induce a splitting of $\Omega^*M
/\ker d$ in $[E] \oplus [F]$ (possibly degenerated) with $[E] = \im
[\Pi_E] = (E + \ker d) /\ker d$ and $[F]= \im [\Pi_F] = (F+ \ker d)
/\ker d$.
\begin{prop}
  The cut-off cones $C_\varepsilon (E(\Delta \leq 1),d) \subset
  \Omega^*M /\ker d \simeq (\ker d)^\bot$ are uniformly shrinking
  around $[E]$ relatively to $[F]$. Precisely, $\exists C $ such that
  $$\|\Pi_F \alpha\|_{\Omega^*M /\ker d } \leq C \varepsilon
  \|\alpha\|$$
  for all $\alpha \in C_\varepsilon (E(\Delta \leq
  1),d)$.
\end{prop}
\begin{proof}
  Recall that $\Pi_F = Qd+dQ$, so that for $\alpha \in C_\varepsilon
  (E(\Delta \leq 1),d)$ one has
  \begin{align*}
    \|\Pi_F \alpha\|_{\Omega^*M /\ker d } &\leq \|Qd \alpha\|\\
    & \leq C \|d\alpha\| \leq C\varepsilon \|\alpha\|,
  \end{align*}
  where $C$ is a bound for $Q$ on $E(\Delta \leq 1) \supset d(E(\Delta
  \leq 1))$.
\end{proof}
\begin{rems}
\item - Notice that $[\Pi_E]$ and $[\Pi_F]$ are not orthogonal
  projections.
\item - As $F= \im d_0^{-1} + \im dd_0^{-1}$ is pinched between $\im
  d_0^{-1} $ and $\im d_0^{-1} + \im d$ one has $[F] = [\im
  d_0^{-1}]$. Also $E= \ker d_0^{-1}\cap \ker d_0^{-1} d \subset \ker
  d_0^{-1}$ so that $[E] \subset [\ker d_0^{-1}]$.
\end{rems}

We close this section with an application of theorem \ref{thm:3:3} to
$d_c$-harmonic decomposition on \emph{compact} C-C manifolds.
\begin{prop}\label{prop:3:6}
  Let $M$ be a compact $E_0$-regular C-C manifold. Then $\im d_c$ and
  $\im \delta_c$ are closed in $L^2$ and we have an orthogonal
  Hodge-de Rham splitting
  \begin{displaymath}
    E_0 = \mathcal{H}_c \oplus \im d_c \oplus \im \delta_c,
  \end{displaymath}
  where $\mathcal{H}_c = \ker d_c \cap \ker \delta_c$, is the space of
  $d_c$-harmonic forms, isomorphic to de Rham's cohomology of $M$.
\end{prop}
\begin{rem}\label{rem:3:7}
  This $d_c$-harmonic representation of the cohomology is probably
  more interesting when $E_0$ is of pure C-C weight (like on one
  forms), since in that case $\mathcal{H}_c$ will be invariant under
  the natural C-C dilations of the riemannian metric $g\rightarrow
  \lambda g_{V_1} + \lambda^2 g_{V_2} + \cdots + \lambda^r g_{V_r}$,
  with $V_i = H_i \ominus H_{i-1}$. In that pure weight case, one can
  also obtain analytic information on the regularizing properties of
  $d_c$, showing in particular the smoothness of $\mathcal{H}_c$ (see
  section \ref{sec:5:1}).
\end{rem}
\begin{proof}
  On a compact manifold, $0$ is isolated in the (discrete) spectrum of
  the de Rham's Laplacian $\Delta$. Therefore, on $(\ker d)^\bot$
  $$\|d\alpha\|^2 = (\Delta \alpha, \alpha) \geq \lambda_1
  \|\alpha\|^2$$
  for some $\lambda_1 > 0$. Then $C_\varepsilon(d) =
  \{0\}$ for $\varepsilon < \lambda_1$, and the near-cohomology of $M$
  vanishes. By theorem \ref{thm:3:3} the same is true at the level of
  $(E_0,d_c)$. This implies the same kind of bound for
  $q_c=\|d_c\|^2$, showing that $\im d_c$ is closed, and the same for
  $\im \delta_c$ by $*$-duality.
\end{proof}

\subsection{From cones to spectral densities}
\label{sec:3:2}

An interesting feature of near-cohomology lies in its stability under
general bounded maps, whereas the spectral decomposition of an
operator (like the Laplacian) is only preserved by isometries. This is
the advantage of working with the quadratic form $q(\alpha) =
\|d\alpha\|^2$ on $\Omega^*M /\ker d $, up to multiplicative
constants, instead of doing the (much more) fine study of the exact
spectral decomposition of $\Delta$ on forms.  Yet, the cones
$C_\varepsilon = \{\alpha \mid q(\alpha) \leq \varepsilon^2
\|\alpha\|^2 \}$ encoded the ``rough'' information we are looking for
of the asymptotic spectral density of $\Delta$ at $0$.

We first observe that $E(0< \delta d \leq \varepsilon^2)$ is a closed
linear space in $C_\varepsilon \subset (\ker d)^\bot$. Moreover, one
has $q(\alpha) > \varepsilon^2 \|\alpha\|^2$ on its orthogonal
$E(\delta d > \varepsilon^2) = E(0< \delta d \leq \varepsilon^2)^\bot
\cap (\ker d)^\bot$.  Therefore, $C_\varepsilon \cap E(\delta d >
\varepsilon^2) = \{0\}$. A consequence is that \emph{any} linear space
$L$ lying in $C_\varepsilon$ projects \emph{injectively} into $E(0<
\delta d \leq \varepsilon^2)$, showing in some sense that this space
is the ``largest'' possible one inside $C_\varepsilon$.

\subsubsection{$\Gamma$-dimension and the rational case}

In order to translate numerically this observation, one needs a notion
of dimension working on some infinite dimensional linear spaces.
%, such that $\dim
%L \leq \dim L'$ if $L \subset L'$ and $\dim f(L) = \dim L$ if $f$ is
%bounded and injective. 
There is an classical one, called $\Gamma$-dimension, in the case a
discrete group of isometries $\Gamma$ is acting co-compactly on the
manifold (eg a Galois covering of a compact manifold). A complete
exposition of this may be found in Atiyah's original work \cite{At}
(or \cite{Pa} for a survey of its properties and applications).

Briefly, in the case we are concerned in, let $M$ be a compact
manifold and $\tilde M$ its universal cover. Let $P$ be a positive
$\Gamma$-invariant Hermitian operator acting on some
$\Gamma$-invariant Hilbert space $H$ over $\tilde M$. We assume $H$ to
be a subspace of $L^2(\tilde M, V)$, where $V$ is the pull-back of
some finite dimensional vector bundle on $M$. Then $H$ admits a
$\Gamma$-invariant Hilbert base $(e_i)_{i\in I}$ on which
$\Gamma$-acts freely ($H$ is called a free $\Gamma$-module). Choose
one vector by $\Gamma$-orbit. (For instance one can take an Hilbert
base of $H$ restricted to a fundamental domain $\mathcal{F}$ of
$\tilde M$.) The $\Gamma$-trace of $P$ is defined by
\begin{displaymath}
  \tr_\Gamma (P) = \sum_{i\in \Gamma \backslash I} (Pe_i,e_i),
\end{displaymath}
The definition is shown to be independent of the choices of $e_i$
(lemma 2.2 in \cite{Pa}).  If $L$ is a closed $\Gamma$-invariant
subspace of $H$, its $\Gamma$-dimension is defined by $\dim_\Gamma L=
\tr_\Gamma \Pi_L$, where $\Pi_L$ denotes the orthogonal projection
onto $L$.

The basic properties of $\dim_\Gamma$ we need here are
\begin{itemize}
\item first its invariance under $\Gamma$-invariant bounded injection
  (proved using polar decomposition see eg lemma 2.3 in \cite{Pa}),
\item its link with kernel density in the case of a smoothing
  projection $\Pi_L$. For example, if $L=E(\Delta \leq \lambda)$ then
  \begin{equation}\label{eq:3:2:1}
    \dim_\Gamma L = \tr_\Gamma \Pi_L = \int_\mathcal{F} \tr (K(x,x))
    \mathrm{dvol},
  \end{equation}
  where $K(x,y)$ is the smooth Schwartz kernel of $\Pi_L$, and $\tr$
  is the standard finite dimensional trace on $\mathrm{End}(\Lambda^*
  T^*M)$.
\end{itemize}

Applying this to a previous remark that any linear subspace in
$C_\lambda$ projects injectively into $E(0 <\delta d\leq \lambda^2)$,
one obtain the following variational principle (\cite{Gr-Sh})
\begin{equation}
  \label{eq:3:2:2}
  F_{\delta d
    } (\lambda^2) = \dim_\Gamma E(0 <\delta d\leq
  \lambda^2) = \sup_{L \in \mathcal{L}_\lambda}\dim_\Gamma L,
\end{equation}
where $\mathcal{L}_\lambda$ is the set of all $\Gamma$-invariant
closed linear subspaces $L\subset C_\lambda(d)$. Hence the spectral
distribution function of $\delta d$ is encoded in the $\Gamma$-linear
``thickness'' of the near-cohomology cones. Moreover, theorem
\ref{thm:3:2} and the invariance of $\dim_\Gamma$ trough injective
maps, implies Gromov-Shubin's result that $\Gamma$-homotopical
equivalent complexes $(E,d)$ and $(E',d')$ on $\tilde M$ must have
equivalent density functions
\begin{equation}
  \label{eq:3:2:3}
  F_d(C^{-1}\lambda) \leq F_{d'}(\lambda) \leq F_d (C\lambda),
\end{equation}
for small $\lambda$ and some $C>0$, where $F_d(\lambda)$ is defined in
general by
\begin{equation}
  \label{eq:3:2:4}
  F_d (\lambda) = \sup_{L \in \mathcal{L}_\lambda}\dim_\Gamma L,
\end{equation}
with the same $\mathcal{L}_\lambda$ as above. In particular they have
the same Novikov-Shubin exponents
\begin{equation}
  \label{eq:3:2:5}
  \alpha_d =  \liminf_{\lambda \rightarrow
  0} \frac{\ln F_d ( \lambda)}{\ln \lambda}
\end{equation}
(but also same limsup, or any other dilatationally invariant limit).
This applies to theorem \ref{thm:3:3}, \emph{if} the manifold admits a
cocompact discrete action.
\begin{defn}
  A graded nilpotent Lie group $G$ is called rational if it admits a
  cocompact discrete group $\Gamma$.
\end{defn}
For such a group rationality is equivalent to being able to find a
basis of $\mathfrak{g}$ with brackets given by rational coefficients
(see \cite{Co-Gr}). Then a cocompact group $\Gamma$ is given by the
exponential of a corresponding integral lattice in $\mathfrak{g}$.
\begin{thm}\cite{Ru}\label{thm:3:9}
  Let $G$ be a rational graded nilpotent Lie group. Then de Rham's
  complex and $(E_0,d_c)$ have the same Novikov-Shubin exponents.
\end{thm}
These exponents do not depend on the choice of $\Gamma$ here.  One can
endow $G$ with a (left) invariant metric, in which case all operators
and spaces become $G$-invariant and formulae \eqref{eq:3:2:2} and
\eqref{eq:3:2:1} reduce to
\begin{displaymath}
  F_{\delta d}(\lambda^2) = \mathrm{vol} (\Gamma \backslash G) \tr
  (K_{\Pi_{E(0 < \delta d \leq \lambda^2)}} (e,e)).
\end{displaymath}
Hence $\dsp \liminf_{\lambda \rightarrow 0} \frac{\ln F_{\delta
    d}(\lambda^2)}{\ln \lambda}$ is independent of $\Gamma$.

We remark also that theorem \ref{thm:3:9} is true on any fundamental
cover of a compact $E_0$ regular C-C manifold. We don't state it in
this generality since for applications we will need some dilation
respecting the C-C structure, that is more or less to work on a graded
group.

\subsubsection{Extension to non-rational groups}

Although ``geometric'' nilpotent groups (associated to semi-simple
geometries) tend to be rational, they are certainly negligible in the
variety of nilpotent Lie groups. Moreover the question of the
asymptotic heat decay makes sense on any group, and it seems unlikely
this rough invariant should depend so sharply on the existence of a
rational structure. (Observe for instance that the value on functions
is the growth rate of $G$, which turns out to be an integer for all
nilpotent Lie groups.) Therefore it is natural to extend the previous
result to non-rational groups. The point is to use a counterpart to
$\Gamma$-dimension for $G$-invariant Hilbert spaces. We describe two
different approaches.

\item $\bullet$ The first one relies on the use of unitary
  representation theory of nilpotent Lie groups.  On such a group $G$,
  the Fourier-Plancherel transform gives an isometry between $L^2(G)$
  with $L^2(\widehat G)$ the unitary dual of $G$ (see eg \cite{Co-Gr}
  for a complete exposition). We briefly explain how the Plancherel
  measure on $\widehat G$ gives rise to a relevant measure on
  $G$-invariant Hilbert spaces in $L^2(G)$.
  
  A bounded (left) $G$-invariant map $P$ acting on $L^2(G)$ has a
  Schwarz kernel $K(x,y)$ (a distribution on $G\times G$) which is
  invariant $K(gx,gy) =K(x,y)$.  Therefore $P$ acts as a convolution
  product $Pf =k * f$ which ``diagonalizes'' trough Fourier-Plancherel
  transform in $\pi (Pf) = \pi (k) \circ \pi (f) $, where $\pi$ is an
  irreducible unitary representation of $G$ (in bijection with
  covectors $\xi$ in $\mathfrak{ g}^*$ modulo the coadjoint action of
  $\mathfrak {g}$ by Kirillov's orbit method).  Recall that for $f\in
  L^2(G)$, $\widehat {f} : \pi \rightarrow \pi (f)$ consists in a
  measurable field of Hilbert-Schmidt operators acting for each $\pi$
  on the associated Hilbert space $H_\pi$. For $f\in L^1(G)$, one has
  \begin{displaymath}
    \pi(f) = \int_G f(x) \pi(x^{-1}) dx.
  \end{displaymath}
  Plancherel formula is
  \begin{equation}\label{eq:plancherel}
    \|f\|_2^2 = \int_G | f(x)|^2 dx =\|\widehat{f}\|_2^2 =
    \int_{\widehat G} \|\pi (f)\|_{HS}^2 d\mu (\pi)
  \end{equation}
  with $\|\pi (f)\|_{HS}^2 = \tr (\pi (f)^* \pi (f) )$, while Fourier
  inversion formula reads
  \begin{equation}\label{eq:fourier}
    f(e) =\int_{\widehat G} \tr (\pi (f)) d\mu (\pi)
  \end{equation} 
  for $f\in \mathcal{S}(G)$. Then to a bounded $P$ is associated a
  field of bounded $\pi (k)\in \mathrm{End} (H_\pi)$.  If moreover $P$
  is an orthogonal projection onto a $G$-invariant closed space $L$,
  these $\pi (k)$ have to be orthogonal projections onto spaces $L_\pi
  \subset H_\pi$. We can now define the $G$-dimension of $L$ as
  \begin{displaymath}
    \mathrm{dim}_G L = \int_{\widehat G} \mathrm{dim} (L_\pi) d\mu (\pi).
  \end{displaymath}
  For a smoothing $\Pi_L$ Fourier inversion formula gives
  \begin{displaymath}
    \mathrm{dim}_G L = k(e) = K_{\Pi_L}(e,e),
  \end{displaymath}
  to be compared with its $\Gamma$-dimension if $G$ is rational. Using
  the previous discussion one find also that $\mathrm{dim}_G$ is
  preserved trough injective bounded $G$-invariant maps has needed.
  All this leads to an extension of theorem \ref{thm:3:9} to any
  graded nilpotent Lie group.
  \begin{thm}\label{thm:3:10}
    Theorem \ref{thm:3:9} holds without rationality assumption,
    replacing the $\Gamma$ by the $G$-dimension.
  \end{thm}

\item $\bullet$ The second approach of the previous result uses more
  elementary tools, and follows Atiyah's presentation of the
  $\Gamma$-trace in \cite{At} section 4. Let $G$ be a graded nilpotent
  group. The main fact we need is that although $G$ has no compact
  quotient in general, it admits a nice covering. This has been
  observed in a close form by Helffer and Nourrigat in \cite{He-No},
  lemma 4.5.2.
  \begin{lemma}\label{lem:3:3}
    Let $G$ be a graded nilpotent Lie group. There are a relatively
    compact domain $D$ and a discrete set $Z \subset G$, such that
    $Z^{-1} =Z$ and
    \begin{displaymath}
      G = \bigsqcup_{z \in Z} zD.
    \end{displaymath}
  \end{lemma}
  \begin{proof}
    The Lie algebra $\mathfrak{g}$ splits in $\mathfrak{g_1} \oplus
    \mathfrak{g_2} \cdots \oplus \mathfrak{g_r}$ with
    $[\mathfrak{g_1}, \mathfrak{g_i}] \subset \mathfrak{g_{i+1}}$.
    Choose an integral \emph{additive} lattice $Z_i \subset
    \mathfrak{g_i}$, together with a fundamental domain $D_i$ of $Z_i$
    in $\mathfrak{g_i}$.  Consider $D= \exp (D_1 \times D_2 \cdots
    \times D_r)$ and $Z=\exp (Z_1 \times Z_2 \cdots \times Z_r)$. Now
    at the Lie algebra level, given $y= (y_1, \cdots, y_r)\in
    \mathfrak{g}$, the equation
    \begin{displaymath}
      (z_1, \cdots, z_r).(x_1,\cdots , x_r) = (y_1, \cdots, y_r)
    \end{displaymath}
    leads to a triangular system whose diagonal terms are
    \begin{displaymath}
      z_i + x_i = y_i + \mathrm{(non\ linear)\ function\ of\  }
      z_j,x_k\  
      \mathrm{with\ }j, k <i. 
    \end{displaymath}
    Hence it admits a unique solution with $z\in Z$ and $x\in D$.
    (This proof adapts on any nilpotent Lie group.)
  \end{proof}
  This gives a discretization of $L^2(G)$ in $\ell^2 (Z) \otimes
  L^2(D)$, which allows to follow the classical construction of
  $\Gamma$-trace as given in \cite{At} section 4 (or \cite{Pa}). For
  instance a bounded $G$-invariant operator $P$ may be called
  $Z$-Hilbert Schmidt if its kernel belongs to $L^2(Z\backslash
  G\times G)$, which identifies either to $L^2(D \times G)$ or
  $L^2(G\times D)$ by the symmetry $Z^{-1} = Z$. Its $Z$-Hilbert
  Schmidt norm may be defined by
  \begin{align*}
    \|P\|_{ZHS}^2 & = \int_{G \times D} |K(x,y)|^2 dx dy =
    \|P\chi_D\|_{HS}^2 = \sum_{e_i}
    \|(P\chi_D) e_i\|^2 \\
    & = \int_{D\times G} |K(x,y)|^2 dx dy = \|\chi_D P\|_{HS}^2 =
    \sum_{e_i} \|(\chi_D P) e_i\|^2 ,
  \end{align*}
  for $e_i$ an Hilbert base of $L^2(G)$ and $\chi_D$ the
  characteristic function of $D$.  It is easily seen to be independent
  on the choice of $D$ and $*$-invariant $\|P\|_{ZHS} =
  \|P^*\|_{ZHS}$. We check the first point, if $D'$ is another
  ``fundamental domain'', one has
  \begin{align*}
    \|P\|_{DHS}^2 & = \int_G \left(\int_D |K(x,y)|^2 \bigl(\sum_{z \in
        Z}
      \chi_{zD'} \bigr) (y)  dy \right) dx \\
    & = \int_G \sum_{z \in Z^{-1}=Z} \int_{zD} |K(x,y)|^2 \chi_{D'}
    (y) dy dx\\
    & = \int_G \bigl(\int_G |K(x,y) |^2 \chi_{D'} (y) dy \bigr)dx =
    \|P\|_{D'HS}^2
  \end{align*}
  A positive hermitian $P$ is of $Z$-trace class if $\tr_Z (P) =
  \sum_{i\in I} (Pe_i, e_i)$ converges where $(e_i)$ an Hilbert base
  of $L^2(D)$, equivalently $\chi_D P \chi_D$ is of trace class on
  $L^2(G)$. The relation with the $Z$-HS class being
  \begin{displaymath}
    \tr_Z (P^*P) = \tr (\chi_D P^*P \chi_D) = \|P\chi_D\|_{HS}^2 =
    \|P\|_{ZHS}^2. 
  \end{displaymath}
  Again if $\Pi_L$ is a smoothing $G$-invariant projection on $L$,
  $\dim_Z L= \tr_Z \Pi_L$ reduces to $\mathrm{vol} (D) K_{\Pi_L}
  (e,e)$ as in the case of $\Gamma$-dimension in \eqref{eq:3:2:1}. So
  we get a version of the previous theorem \ref{thm:3:10} using this
  ``$Z$-dimension'' instead.
  
  Observe that in all considered $\Gamma$, $G$, or $Z$ approaches of
  dim, one can uses only $G$-invariant closed spaces to determine the
  cones thickness function $F(\lambda)$ in formula \eqref{eq:3:2:4}.
  Indeed the maximal (through injective projection) Hilbert space in
  $C_\lambda(d)$ is $E(0<\delta d \leq \lambda^2)$ a $G$-invariant
  space, not only a $\Gamma$ or $Z$-invariant one.

\subsection{The algebraic pinching of heat decay}
\label{sec:3:3}

We can now give a first application of these constructions to our
initial problem of estimating the heat decay on ($1$-)forms of graded
groups. We define the C-C weight of $G$, with respect to its grading,
as $\dsp N(G)=\sum_{k=1}^r k \mathrm{dim} (\mathfrak{g}_k)$. This
$N(G)$ is the weight of $\Lambda^{\mathrm{max}}T^*G$, occurring in the
Jacobian of the dilations of $G$. When the group $G$ is filtered,
$N(G)$ is also its Hausdorff dimension, that is the growing rate of
the volume of large balls in $G$. We also consider the splitting of
$E_0^k= H^k(\mathfrak{g})$, the $k$-forms of $E_0$, in its various C-C
weights $p$
\begin{displaymath}
  E_0^k = \bigoplus_{p=N_k^{\min}}^{N_k^{\max}} E_0^{k,p}.
\end{displaymath}
$E_0^k$ is said of pure weight $p$ if $p= N_k^{\min}= N_k^{\max}$.
\begin{defn}
  We will use $\dsp \beta_k = N(G) / \alpha_k $ as a convenient
  renormalization of $\alpha_k$, the $k$-th Novikov-Shubin exponent of
  $G$ (twice the heat decay on coclosed $k$-forms of $G$ by
  \eqref{eq:1:1}).
\end{defn}
\begin{thm}\cite{Ru}
  \label{thm:3:3:1}(Algebraic pinching of $\beta_k$.) Let $G$ be a
  graded nilpotent Lie group $G$. If $E_0^k= H^k(\mathfrak{g})$ has
  pure weight $N_k$ then,
  \begin{displaymath}
    \delta N_k^\mathrm{min} = \max (N_{k+1}^{\min} - N_k , 1) \leq
    \beta_k \leq N_{k+1}^{\max} - N_k = \delta N_k^\mathrm{max}, 
  \end{displaymath}
  giving an estimation of $\beta_k = N(G) / \alpha_k$ from the weight
  gaps $\delta N_k$.
\end{thm}
\begin{rem}
  This is of course consistent with the case of functions $k=0$
  mentioned in theorem \ref{thm:1:2}. Namely for a filtered Lie group
  (bracket generated by $D=\mathfrak{g}_1$), one has $E_0^1\simeq
  \Lambda^1 D^*$ of pure weight $1$ (section \ref{sec:2:4}). Therefore
  $\beta_0= 1$ and we recover that $\alpha_0= N(G)$ is the growth of
  $G$.
\end{rem}
\begin{proof}
  By theorem \ref{thm:3:10}, we can use $(E_0, d_c)$ to estimate
  $\beta_k$. Then the result relies on an homogeneity argument of
  $d_c$ trough the dilations $h_r$ of $G$. Observe that on $E_0^k$,
  \begin{equation}
    \label{eq:3:3:1}
    d_c = d_c^{\delta N_k^\mathrm{min}} + \cdots + d_c^{\delta
    N_k^\mathrm{max}},
  \end{equation}
  where $d_c^i$ increases the weight by $i$.
  
  We recall that for a form $\alpha$ of weight $N(\alpha)$, the
  integral $L^2$-norms transforms as $\|h^*_r \alpha\|_G= r^{N(\alpha)
    -N(G)/2} \|\alpha\|_G$. We now see how the near-cohomology cones
  $C_\lambda(d_c)$ rescale through $h_r^*$. Firstly, by the
  homogeneity assumption of $E_0^k$, the space $E_0^k/ \ker d_c \simeq
  (\ker d_c)^\bot$ is preserved. By the previous remarks, one has then
  for $\alpha\in E_0^k \cap C_\lambda (d_c) $ and $r\geq 1$,
  \begin{align*}
    \|d_c (h_r^* \alpha)\| = \|h_r^* d_c \alpha\| & \leq
    r^{N_{k+1}^{\max} -N(G)/2} \|d_c \alpha\|\\
    & \leq \lambda  r^{N_{k+1}^{\max} -N(G) /2} \|\alpha\|\\
    & = \lambda r^{\delta N_k^\mathrm{max}} \|h_r^* \alpha\|.
  \end{align*}
  Therefore
  \begin{displaymath}
    h_r^* (C_\lambda(d_c)) \subset C_{\lambda r^{\delta
    N_k^\mathrm{max}}}(d_c).
  \end{displaymath}
  In the opposite direction, for $\varepsilon \leq 1$ and $\alpha \in
  E_0^k \cap C_\lambda (d_c) $
  \begin{align*}
    \| d_c (h_\varepsilon^* \alpha)\| = \|h_\varepsilon^* d_c \alpha\|
    & \leq \varepsilon^{\max (N_k^\mathrm{min}, N_k + 1) -N(G) /2 }
    \|d_c \alpha\| \\
    & \leq \lambda \varepsilon^{\max (N_k^\mathrm{min}, N_k + 1) -N(G)
      /2 } \|\alpha\|\\
    & \leq \lambda \varepsilon^{\delta N_k^\mathrm{min}}
    \|h_\varepsilon^* \alpha\|,
  \end{align*}
  giving
  \begin{displaymath}
    h_\varepsilon^* (C_\lambda(d_c)) \subset C_{\lambda
    \varepsilon^{\delta N_k^\mathrm{min}}} (d_c).
  \end{displaymath}
  Putting this together we have found that, for $r \geq 1$
  \begin{equation}
    \label{eq:3:3:2}
    h_r^* (C_{r^{-\delta
    N_k^\mathrm{max}}}(d_c)) \subset C_1 (d_c) \subset h_r^*
    (C_{r^{-\delta N_k^\mathrm{min} }}( d_c)).
  \end{equation}
  
  The next point to check is the effect of dilations on dimension of
  $G$-invariant Hilbert spaces $L\subset E_0$. We have for $n\in \N$,
  \begin{equation}
    \label{eq:3:3:3}
    \mathrm{dim}_Z (h_n^* L) = \mathrm{dim}_{h_n Z} (L) = n^{N(G)}
    \mathrm{dim}_Z (L). 
  \end{equation}
  We use the discrete dilations $h_n$ because we will need that $h_n
  (Z) \subset Z$. (We could avoid this working with the less
  elementary $G$-dimension).  For $\Gamma$-dimension this formula is a
  direct consequence of its behavior through finite covering. The
  proof is the same here. Since $n^{N(G)/2 - N_k } h_n^* $ is an
  isometry on $L^2(E_0^k)$, we have
  \begin{displaymath}
    \Pi_{h_n^* L} = h_n^* \Pi_L h_n^{*-1}.
  \end{displaymath}
  Therefore,
  \begin{align*}
    \mathrm{dim}_Z (h_n^* L) &= \tr_Z (\Pi_{h_n^* L}) = \tr (\chi_D
    h_n^* \Pi_L h_n^{*-1} \chi_D) \\
    &= \tr (h_n^* \chi_{h_n D} \Pi_L \chi_{h_n D} h_n^{*-1}) = \tr (
    \chi_{h_n D} \Pi_L \chi_{h_n D} ) \\
    & = \mathrm{dim}_{h_n Z} (L),
  \end{align*}
  Now by lemma \ref{lem:3:3}, $h_n D$ is a ``fundamental domain'' for
  the new ``integral lattice'' $h_n Z$ which splits, for $n\in \N$, in
  $n^{N(G)}$ domains $D'$ for $Z$, giving \eqref{eq:3:3:3}. This is
  applied to $G$-invariant Hilbert linear sub-spaces $L\subset
  C_\lambda (d_c)$ and leads with \eqref{eq:3:3:2} and
  \eqref{eq:3:2:4} to the following pinching for the cone thickness
  function $F_{d_c}$
  \begin{displaymath}
    F_{d_c} (n^{-\delta N_k^\mathrm{max}}) \leq n^{-N(G)} F_{d_c}(1) \leq
    F_{d_c} (n^{-\delta N_k^\mathrm{min}}),
  \end{displaymath}
  for $n\in \N$ and finally
  \begin{equation}
    \label{eq:3:3:4}
    C F_{d_c}(1) \lambda^{N(G)/ \delta N_k^\mathrm{min}} \leq
    F_{d_c}(\lambda) 
    \leq C' F_{d_c}(1) \lambda^{N(G) /\delta N_k^\mathrm{max}},
  \end{equation}
  for $\lambda > 0$ (even real) and some positive constants $C,C'$.
  
  The last observation is that all these numbers are finite and
  strictly positive.  Finiteness comes from the fact that this is true
  ``upstairs'' for the homotopically equivalent de Rham complex.
  Indeed by ellipticity of the de Rham Laplacian, spectral projectors
  have smooth kernels, and therefore of finite diagonal trace.
  Non-vanishing amounts to the fact that $0$ belongs to the spectrum
  of $\Delta$ on $k$-forms for any $k$ on nilpotent groups. This can
  be viewed at the $(E_0, d_c)$ level. The cones $C_\lambda(d_c)$
  don't vanish, because any sufficiently dilated non-closed smooth
  form of $E_0$ lies in it, as $d_c$ (unlike $d$) strictly increases
  the weight. We note that $d_c$ can't vanish identically in degree $k
  < \mathrm{dim}G$, because this complex is a resolution of functions
  of $G$ and has to be of maximal length.
\end{proof}

We can give a more geometric meaning of the previous theorem
\ref{thm:3:3:1} in the case of one forms. Namely we have seen in
section \ref{sec:2:4} that when the group is filtered,
\begin{itemize}
\item firstly $E_0^1$ identifies with $\Lambda^1 D^*$ forms of weight
  $1$,
\item $E_0^2 = H^2(\mathfrak{g})$ may be interpreted as $\Lambda^1
  r^*$, the dual to the space $R= N / (N, \widetilde G)$ of generators
  of the relations $N$ of $G$, with respect to the free Lie group
  $\widetilde G$.
\end{itemize}
Therefore theorem \ref{thm:3:3:1} translates into the following
result.
\begin{cor}\label{cor:3:3:2}
  Let $G$ be a filtered nilpotent Lie group. Then $\beta_1 + 1 = N(G)
  / \alpha_1 + 1$ is pinched between the minimal and maximal weight of
  the relations giving a presentation of $G$ from the free Lie group.
  In particular $1\leq \beta_1 \leq r$ if $G$ has rank $r$.
\end{cor}

\section{Examples}
\label{sec:3:4}

The previous result gives an estimation of heat decay from the
knowledge of the weights of the Lie algebra cohomology
$H^*(\mathfrak{g})$. This pinching relies therefore on finite
dimensional linear systems computations (that may of course become
incomprehensible as the dimension increases).

\subsection{Quadratically presented groups}
\label{sec:4:1}

The first class of nilpotent groups on which the previous results
applies sharply are quadratically presented groups. At the Lie algebra
level they are groups admitting a presentation using only linear
combinations of bracket of first order in the generators
$D=\mathfrak{g}_1$. By the previous discussion this may also be
checked dually by showing that $H^2(\mathfrak{g})$ contains only
$2$-forms of weight $2$ (horizontal $2$-forms in $\Lambda^2 D^*$).

For all these groups one has $\beta_1= 1$, meaning that $\alpha_1 =
N(G)$ (in the strong sense that on $\Omega^1 G$, the heat kernel
satisfies $K_t(x,x) \asymp t^{-N(G)/2}$ when $t\rightarrow +\infty$,
as follows from \eqref{eq:3:3:4}). Thus heat decay on one forms is the
quickest possible and the same as on functions.  Fortunately these
groups exist.

\item $\bullet$ The simplest examples are given by the Heisenberg
  groups $\mathbf{H}^{2n+1}$ of dimension $2n+1\geq 5$.  They may be
  presented by the relations $[X_i, X_j] = [Y_i, Y_j] = 0$ with $1
  \leq i, j\leq n$, $[X_i, Y_j] = 0$ for $i \not=j$ and (\emph{for}
  $n> 1$), $[X_i, Y_i] = [X_j, Y_j]$ ($ =T$ the generator of the
  center).  In fact using the Lie algebra computation in section
  \ref{sec:2:3} and theorem \ref{thm:3:3:1}, we have $\alpha_p = N(G)$
  for all $0 \leq p \leq 2n$ except $\alpha_n= N(G)/2$. In particular
  the value $\alpha_1 (\mathbf{H}^3) =2$ (first computed by Lott
  \cite{Lo}) is half the Hausdorff dimension. This corresponds
  geometrically to the fact that $\mathbf{H}^3$ has a \emph{cubic}
  presentation : $[X,[X,Y]] = [Y, [X,Y]]= 0$.
  
\item $\bullet$ There are in fact many other examples. In \cite{Gr}
  section 4.2, Gromov introduced the class of two step graded groups
  that possess a $\Omega$-regular legendrian two plane $P$. That means
  there exists $X_1, X_2$ in $D$ such that $[X_1,X_2]= 0$ and the map
  \begin{align*}
    \Omega: D & \rightarrow \mathrm{End}(P, [D, D]) \\
    X & \rightarrow [X,\cdot ]\ \mathrm{\ acting\ on\ } P=
    \mathrm{span} (X_1, X_2),
  \end{align*}
  is surjective. (Giving $T_1, T_2 \in \mathfrak{g}_2 = [D,D]$, there
  exists $X$ in $D$ such that $[X, X_1] = T_1$ and $[X, X_2] = T_2$.)
  One can show easily (see \cite{Gr} 4.2 A'') that this condition is
  generically satisfied (is a Zariski open dense set) for graded
  groups associated to distributions $D^k \subset \R^n$ as long as $n
  \leq 3k/2 - 1$. (The necessity of this bound is a dimension counting
  in the previous formula noting that $P \subset \ker \Omega$.)
  
  The geometric interest of this class of groups is that Gromov has
  shown they are quadratically fillable in the sense that closed
  horizontal curves $\gamma$ can be filled by horizontal disks of area
  controlled by $\mathrm{length}(\gamma)^2$. (This is a global
  non-linear hard property relying infinitesimally on the presence of
  a lot of flexible legendrian planes.) Our observation is simply that
  \begin{prop}
    $\Omega$-regular groups are quadratically presented, and thus have
    $\beta_1 =1$.
  \end{prop}
  \begin{proof}
    We have to show that $H^2(\mathfrak{g})$ doesn't contain $2$-forms
    of weight $3$, that is the map
    \begin{equation}\label{eq:3:4:1}
      \begin{aligned}
        d_0 : \Lambda^{2,(3)} \mathfrak{g}^* = D^* \wedge
        \mathfrak{g}_2^*
        & \longrightarrow  \Lambda^{3,(3)} \mathfrak{g}^* =\Lambda^3 D^* \\
        \sum \alpha_i \wedge \theta_i & \longrightarrow - \sum
        \alpha_i \wedge d_0 \theta_i = \sum \alpha_i \wedge \theta_i
        ([\cdot, \cdot])
      \end{aligned}
    \end{equation}
    is injective. Let $\beta = \sum \alpha_i \wedge \theta_i$ be in
    $\ker d_0$.  We first show that $\alpha_i$ vanishes on $P=
    \mathrm{span} (X_1, X_2)$. We fix a base $\theta_i$ of
    $\mathfrak{g}_2^*$ in the previous formula. Given $T_{i_0}$ dual
    to $\theta_{i_0}$ we can find by $\Omega$-regularity a vector
    $X\in D$ such that $[X, X_1]= 0$ and $[X,X_2] = T_{i_0}$ (and
    still $[X_1,X_2] = 0$). Hence
    \begin{displaymath}
      0  = d_0 \beta (X_1, X_2, X)  = \sum \alpha_i (X_1)
      \wedge \theta_i ([X_2, X]) = - \alpha_{i_0} (X_1), 
    \end{displaymath}
    and similarly on $X_2$. Now taking $X$ as before and $Y\in D$ such
    that $[Y, X_2] =0$, one gets
    \begin{displaymath}
      0  = d_0 \beta (X_2, X, Y)  = \sum \alpha_i (Y) \wedge \theta_i
      ([X_2, X])  = - \alpha_{i_0} (Y), 
    \end{displaymath}
    so that $\alpha_{i_0}$ vanishes on $\ker\ [\ \cdot\ , X_2]$. For
    the last step, choose $Y\in D$ such that $[Y, X_1] = T_{i_0}$ and
    $[Y, X_2]= 0$ then, for any $X\in D$
    \begin{displaymath}
      0 = d_0 \beta  (X_1, X, Y) = \sum \alpha_i (X)
      \wedge \theta_i ([Y, X_1]) = \alpha_{i_0} (X).
    \end{displaymath}
    This shows that $\alpha_{i_0} = 0$, and finally $\beta = 0$.
  \end{proof}
  
  The previous fact may be seen more geometrically using Gromov's
  result that these groups are quadratically fillable. If any
  horizontal curve can be filled by an horizontal disk then certainly,
  for any $1$-form $\alpha$ on $G$, the vanishing of $d\alpha$ on
  $\Lambda^2 D$ is sufficient to conclude that $\alpha$ is exact.
  Namely $f(m) = \int_\gamma \alpha$ where $\gamma$ is any horizontal
  path from $0$ to $m$, is well defined in that case and satisfies $d
  f = \alpha$ along $D$. Hence $\theta = \alpha - df$ is a vertical
  form such that $d_0 \theta = d \alpha$ on $\Lambda^2D$ vanishes, and
  finally $\alpha =df$. This argument means that the quadratic filling
  property implies that they are no non zero closed two forms of
  weight $>2$.  In particular $E_0^2 = H^2(\mathfrak{g})$ has to be of
  pure weight $2$, because invariant $2$-forms on $G$ coming from $H^2
  (\mathfrak{g}) = \ker d_0 / \im d_0$ are always closed. We recover
  therefore that $G$ is quadratically presented.

\item $\bullet$ Gromov's quadratically fillable groups have ``a lot''
  of integrable legendrian planes $P= (X, Y) \subset D$ which lead to
  simple quadratic relations $[X, Y] = 0$ and allows a quadratic (in
  weight) integration of $2$-forms. Anyway, the presence of these flat
  planes is not a necessary condition for a group to have a quadratic
  presentation. Here is a ``sporadic'' example with no flat plane at
  all but which is quadratically presentable anyway. This is a
  $J$-group corresponding to a $D^8 \subset \R^{15}$. We have to
  describe a $7$-dimensional subspace $L$ of curvature in $\Lambda^2
  \R^8$, such that
  \begin{displaymath}
    \theta \in \R^7 \rightarrow d\theta = g_D(J(\theta)\  ,\  ) \in L  
  \end{displaymath}
  with $- J(\theta)^2 = \|\theta\|^2 \id = J(\theta)^* J(\theta)$.
  This map $J$ can be realized either as a particular Clifford action,
  or the left multiplication of the imaginary Cayley numbers $\im
  \mathbb{O}$ on $\mathbb{O}$. It can been checked, using direct
  calculations with the $J_i$, $1\leq i\leq 7$ given in tables of
  Clifford algebras (see eg \cite{Hu} section 11), or representation
  theory (see \cite{Sa}) that the map $d_0$ in \eqref{eq:3:4:1} is
  injective (even an isomorphism, since both spaces have dimension
  56). Therefore this group is quadratically presented, although it
  contains no integrable horizontal plane since $[X,Y]=0$ implies that
  $X \bot J_i Y$ for $1\leq i \leq 7$, which gives collinear $X$ and
  $Y$.
  
\item $\bullet$ The previous example $D^8 \subset \R^{15}$ goes beyond
  the generic bound $n \leq 3 k/ 2 -1$ under which a generic $D^k
  \subset \R^n$ is $\Omega$-regular hence quadratically presented.
  Unfortunately, due to the eight periodicity of their structure,
  Clifford modules fail to produce a linearly increasing number of
  orthogonal complex structures on $\R^k$ for $k > 8$ ($\R^{16k}$ has
  only $8$ more such complex structures than $\R^k$). In fact it
  doesn't seem easy to improve the previous bound $n \leq 3 k/2 -1$ in
  general.  Here is another (too rough?) method of construction of
  quadratic groups, independent of $\Omega$-regularity, but which
  leads to the same bound.
  \begin{prop}
    Any increasing family of extensions of $D =\R^k$ by
    $$L_p = \mathrm{span} (\omega_1, \omega_2, \cdots, \omega_p)
    \subset \Lambda^2 D^*,
    $$
    leads to quadratically presented groups if
    \begin{enumerate}
    \item each $\omega_p$ is of rank $\geq 2(p+1)$,
    \item $\omega_p \notin L_{p-1} + \{\alpha \wedge \beta \in
      \Lambda^2 D^*\}$.
    \end{enumerate}
  \end{prop}
  \begin{proof}
    As in \eqref{eq:3:4:1}, we have to show that $\dsp \sum_{i=1}^p
    \alpha_i \wedge \omega_i = 0$ with $\alpha_i \in D^*$ implies
    $\alpha_i = 0$. This is clear for a single term by assumption
    $\omega_1$ of rank $\geq 4$. Suppose we have a non trivial
    relation $\dsp \alpha_p \wedge \omega_p = \sum_{i=1}^{p-1}
    \alpha_i \wedge \omega_i$, with $\alpha_p \not= 0$ by recurrence
    hypothesis.
    
    We will show that all $\alpha_i$ are multiple of $\alpha_p$.
    Firstly $\omega_p \wedge \alpha_1 \wedge \alpha_2 \wedge \cdot
    \wedge \alpha_p =0$, and since $\mathrm{rk} (\omega_p) \geq
    2(p+1)$, this implies $\alpha_1 \wedge \alpha_2 \wedge \cdot
    \wedge \alpha_p =0$. Therefore some $\alpha_{i_0}$ with $1 \leq i
    < p$ is a linear combination of the others. Now $\omega_p
    \bigwedge_{i \not= i_0} \alpha_i =0$ implies again $\bigwedge_{i
      \not= i_0} \alpha_i =0$ and gives another dependent
    $\alpha_{i_1}$. We repeat this until we obtain $\dsp \alpha_p
    \wedge \omega_p = \sum_{i=1}^{p-1} c_i \alpha_p \wedge \omega_i$
    for some constants $c_i$. This implies that $\dsp \omega_p =
    \sum_{i=1}^{p-1} c_i \omega^i + \alpha_p \wedge \beta_p$, which is
    contrary to the assumption $2$.
  \end{proof}
  \begin{rems}
    - Notice that the first hypothesis implies, and is generically
    satisfied if $2p+2 \leq k$, while the second one is generic if
    $p-1 + 2(k-2) +1 < \frac{k(k-1)}{2} = \mathrm{dim }(\Lambda^2
    D^*)$.  This is seen to be implied by the first condition, which
    amounts finally to $n= p + k \leq 3k /2 -1$ as with
    $\Omega$-regularity.
    
    - The second condition means that $L_p$ doesn't contain any rank
    $2$ form $\alpha \wedge \beta$.  This is certainly a necessary
    condition for quadratic presentation.  Otherwise $\alpha \wedge
    \theta$ and $\beta \wedge \theta$, with $d_0 \theta = \alpha
    \wedge \beta \in L_p$, will be a $d_0$-closed but not exact
    $2$-forms of weight $3$ in the corresponding extension $G$ of $D$.
    Another equivalent formulation is that $G$ can't admit a
    surjective morphism onto the $3$-dimensional Heisenberg group
    $\mathbf{H}^3$.
  \end{rems}
  
\item $\bullet$ So far we have only considered $2$-step groups, but
  quadratically presented groups of arbitrarily high rank do exist.
  Interest in that field came from its relation with Sullivan's
  rational homotopy theory, where these groups are $1$-formal
  $1$-minimal. In particular since K\"ahler manifolds $M$ have formal
  minimal models, the Malcev's completion $\Pi_1(M)\otimes \Q$ of
  $\Pi_1 (M)$ are examples of such groups, at the condition they are
  finite dimensional, that is nilpotent. We refer to \cite{Gr-Mo} for
  details on these problems.
  
  The simplest example of a rank $3$ quadratically presented group is
  given by a generic $3$-step distribution $D^5 \subset \R^8$ studied
  in \cite{Ca-To}. Its structure is given by $\dsp D=
  \mathrm{span}_{1\leq i \leq 5} (X_i)$, $\mathfrak{g}_2 = [D,D] =
  \mathrm{span} (Y_1, Y_2)$, and $\mathfrak{g}_3 = \R Z$ with the dual
  relations
  \begin{displaymath}
    d y_1 = x_1 \wedge x_3 + x_2 \wedge x_4,\  d y_2 = x_1 \wedge x_4
    + x_2 \wedge x_5, \ d z = x_1 \wedge y_1 + x_2 \wedge y_2.
  \end{displaymath}
  Observe that the third cubic form $x_1 \wedge y_1 + x_2 \wedge y_2$
  is closed from the two first one, and thus may be written $d z$ in
  order to kill it in $H^2(\mathfrak{g})$. It turns out (see
  \cite{Ca-To}) that there are no other weight $\geq 3$ $d_0$-closed
  $2$-forms, and that this group is quadratically presented.

\item $\bullet$ Lastly, a series of arbitrarily high rank
  quadratically presented group has been described by Chen in
  \cite{Ch}. Given $(n,k) \in \N \times \N^*$, these $k+1$-rank groups
  are generated by $D=\{ X_1, \cdots, X_n, Y_\alpha, |\alpha| = k\}$,
  with $\alpha = (\alpha_1, \cdots, \alpha_n)\in \N^n$ and $|\alpha| =
  \alpha_1 + \cdots + \alpha_n$. The group structure is given by
  \begin{displaymath}
    d x_i =0,\quad dy_\alpha = 0\ \mathrm{if}\ |\alpha| =k,\quad
    dy_\alpha = \sum y_{\alpha + e_i} \wedge x_i\ \mathrm{if}\ |\alpha|
    < k,
  \end{displaymath}
  with $e_i =(0, \cdots, 1,\cdots, 0)$. One has $\mathfrak{g}_j =
  D^{[j]} = \mathrm{span} (Y_\alpha, |\alpha| = k-j+1)$. These
  examples have still $\beta_1 =1$, showing that this function may be
  independent of the rank of the group.

\subsection{Some higher rank groups}
\label{sec:4:2}

We now give examples of groups with relations of higher weight.

\item $\bullet$ Let $G^{k,r}$ be the ``free'' nilpotent group of rank
  $r$ with $k$ generators. This is the quotient of the infinite
  dimensional free Lie group $\widetilde{G}^k$ in $k$ generators by
  all its elements of weight $> r$. Is is the largest $r$-step
  filtered group, since any other one may be presented as a quotient
  of it. By definition the relations of $G^{k,r}$ with respect to
  $\widetilde{G}^k$ are generated by all words of weight $r+1$ in
  $\widetilde{G}^k$. Therefore corollary \ref{cor:3:3:2} gives
  $$\beta_1 (G^{k,r})= r,$$
  the upper possible bound for $r$-step
  groups, and thus the lowest possible $\alpha_1 =N(G)/r$ and heat
  decay $N(G)/2r$ (with respect to the Hausdorff dimension $N(G)$).

\item $\bullet$ These $G^{k,r}$ are helpful in precising the pinching
  of $\beta_1$ for groups associated to generic distributions $D^k
  \subset \R^n$.  Let $n(k,r) = \mathrm{dim} (G^{k,r})$ (a general
  formula using M\"obius function may be found eg in \cite{Gr},
  section 4.1.B).
  \begin{prop}\label{prop:4:4}
    For a fixed $k$, one has $r-1 \leq \beta_1 (G) \leq r$ for groups
    $G$ associated to a generic $D^k \subset \R^n$, if $r$ is such
    that $n(k, r-1) \leq n < n(k, r)$.
  \end{prop}
  \begin{proof}
    In fact these groups have only relations of weight $r$ and $r+1$
    generically.
    
    Let $(X_i)_{1\leq i\leq k}$ be a polynomial germ of $k$
    independent vectors fields around $0$ in $\R^n$. Define $D^k =
    \mathrm{span} (X_i, 1\leq i \leq k)$, and $G$ the graded Lie group
    associated to this C-C structure at $0$ (by the quotiented Lie
    bracket $[\ ,\ ]_0$ see section \ref{sec:2:2:1}). Recall that if
    $H_i = [H , H_{i-1}]$, and $H_1 = D^k$ one has $\mathfrak{g} =
    \mathfrak{g}_1 \oplus \cdots \oplus \mathfrak{g}_r$, with
    $\mathfrak{g}_i = H_i /H_{i-1}$, at $0$. Generically, on a Zariski
    open dense set in the $(j-1)$-jet of $(X_i)$ at $0$, one has
    $$\mathfrak {g}^j = \mathfrak{g}_1 \oplus \cdots \oplus
    \mathfrak{g}_j$$
    of maximal possible dimension, that is $\min
    (n(k,j), n)$. Therefore $\mathfrak {g}^j$ is isomorphic for $j
    \leq r-1$ to the Lie algebra of $G^{k,j}$, and $G$ has no relation
    of weight $< r$. For $j=r$, one has $\mathfrak {g}^r$ of dim $n$
    (meaning $D^k$ is an H\"ormander $r$-step distribution in $\R^n$),
    and then $\mathfrak{g} = \mathfrak{g}^r$ is of rank $r$. This
    gives only relations of weight $\leq r+1$.
  \end{proof}

\item $\bullet$ In the opposite direction to the previous ``generic''
  situation, nilpotent groups associated to complex semi-simple
  geometries may have increasing rank but still a presentation of
  bounded (although not quadratic) depth.
  
  For instance, consider $N_n \subset SL(n,\C)$ the nilpotent group of
  $\id + $strictly upper triangular matrices. It is a rank $(n-1)$
  group whose Lie algebra is generated by $D=\mathrm{span} (X_i =
  E_{i,i+1}\in M_n (\C), 1 \leq i\leq n-1)$.  $N_3= \mathbf{H}^3$ is
  the $3$ dimensional Heisenberg group, already encountered (it is
  cubically presented and has $\beta_1=2$. For all $n \geq 4$, $N_n$
  has a presentation with both quadratic an cubic relations, which are
  $$[X_i, X_j] = 0 \ \mathrm{for}\ j-i > 1,\ \mathrm{and}\ [X_i, [X_i,
  X_{i+1}]] = [X_{i+1}, [X_i, X_{i+1}]] = 0.$$
  Dually, $H^2(N_n)$ is
  generated by $\theta_{X_i}\wedge \theta_{X_j}$ for $j-i > 1$ and
  $\theta_{X_i} \wedge \theta_{Y_i}$, $\theta_{X_{i+1}} \wedge
  \theta_{Y_i}$ where $\theta_{Y_i}$ is the dual form to $Y_i = E_{i,
    i+2} = [X_i, X_{i+1}]$. All this comes from the general
  description by Kostant of the structure of the cohomology of
  (maximal) nilpotent Lie algebras in semi-simple Lie algebra (which
  splits in multiplicity one factors through the adjoint action of the
  maximal torus) (see \cite{Ca-Os},\cite{De-Si}).
  
  Therefore, for these groups of increasing rank, we have $1 \leq
  \beta_1 \leq 2$. We see that $\beta_1 (N_4) = 2$.
  \begin{proof}
    For $f\in C^\infty_0 (N_4)$, consider $\alpha = f \theta_{X_2}$.
    We observe that $d_c \alpha $ has weight $3$. This is because the
    weight $2$ part of $E_0^2 = H^2 (\mathfrak{n}_4)$ is spanned by
    $\theta_{X_1} \wedge \theta_{X_3}$, and thus
    $$(d_c \alpha)_2 = d\alpha (X_1, X_3) = X_1. \alpha(X_3) -X_3.
    \alpha (X_1) =0.$$
    As a consequence, for a generic $f$, we have
    found a non-closed smooth $L^2$ form in $E_0^1 = \Lambda^1 D^*$
    such that $d_c \alpha$ has weight $3$. We will see in proposition
    \ref{prop:5:5} this implies $\beta_1 \geq 2$.
  \end{proof}
  
\item $\bullet$ We close this section with a first example showing
  that $\beta_1 = N(G) / \alpha_1 $ is not necessarily an integer when
  $N(G)$ is the Hausdorff dimension of a filtered group $G$.
  
  We consider the $4$-dimensional Engel group $G$ (studied in section
  \ref{sec:2:3}). We observe that with its filtered weight $N$, one
  has $N(G) = 2 + 2+ 3 = 7$, while $H^2(\mathfrak{g})$ is generated by
  $\theta_Y \wedge \theta_Z$, of weight $3$ and $\theta_X \wedge
  \theta_T$, of weight $4$. Hence, one has the pinching $2 \leq
  \beta_1^N (G) \leq 3$, leading to
  \begin{equation}
    \label{eq:G1}
    7/3 \leq \alpha_1 (G) \leq 7/2.
  \end{equation}
  
  Now, we can change the grading, using instead $N'$ such that
  $N'(X)=1$, $N'(Y) = 2$, $N'(Z) = 3$ and $N'(T) =4$. In that case
  $H^2 (\mathfrak{g})$ becomes homogeneous of pure weight $5$, while
  $E_0^3 = H^3 (\mathfrak{ g}) = \mathrm{span} (\theta_{X \wedge Z
    \wedge T}, \theta_{Y \wedge Z \wedge T})$ has now mixed weight $8$
  and $9$. Applying theorem \ref{thm:3:3:1} at the level of $2$ forms,
  we find that $8-5 = 3\leq \beta_2^{N'} (G) \leq 9-5 = 4$, giving
  \begin{equation}
    \label{eq:G2}
    10/ 4 \leq \alpha_2(G) \leq 10/3
  \end{equation}
  since $N'(G) =10$. This gives a pinching of the spectral density of
  $\delta d$ on $\Omega^2 (G) / \ker d$ which is $*$-conjugated to $d
  \delta$ acting on $\Omega^2 (G) / \ker \delta$. This last one is
  itself conjugated by $\delta$ to the spectrum of $\delta d$ on
  $\Omega^1 (G) / \ker d$. We get finally that $\alpha_2 (G) =
  \alpha_1 (G)$ and that the second pinching \eqref{eq:G2} is strictly
  sharper than the first one \eqref{eq:G1}. In particular $\beta_1^N
  (G) = 7 / \alpha_1 (G)$ can't be an integer using the ``Hausdorff''
  scaling $N$, the relevant one on functions.

\item $\bullet$ Still on that group, we can use the computations of
  section \ref{sec:2:3} to understand more precisely some aspect of
  the asymptotic heat diffusion on its $1$-forms. We restrict our
  study to forms on $G$ whose components along the left invariant
  vectors fields $X, Y, Z, T$ are invariant functions along $Z,T$,
  that is pull-back functions by $\Pi: G \rightarrow \R^2 = G/ (Z,T)$.
  We note $\Omega_{\R^2}^* G = C(\R^2) \otimes
  \Lambda^*\mathfrak{g}^*$ this class of forms, and look at the
  asymptotic behavior of the ``commutative'' heat, or spectrum of
  $\Delta_G$ as acting on $L^2_{\R^2} (\Omega_{\R^2}^* G)$. Like in
  the true $L^2(G)$ situation, the de Rham complex on $G$ has the same
  restricted $\R^2$-asymptotics as $(E_0, d_c)$ (equivalent
  $\R^2$-near cohomology in fact, as given by the same proofs as
  before).
  
  Now, the action of $d_c$ on $\Omega_{\R^2}^1 G$ may be read from
  \eqref{eq:2:2:6}, by cutting the $Z$ and $T$ arrows. We see that
  $d_c$ factorizes here in $Pd_{\R^2}$ where
  \begin{displaymath}
    d_{\R^2} (f \theta_X + g 
    \theta_Y) = (X.g -Y.f) \theta_Z, 
  \end{displaymath}
  and
  \begin{displaymath}
    P (f\theta_Z)  = (X^2.f) \theta_{X\wedge T} + (Y.f)
    \theta_{Y\wedge Z}.
  \end{displaymath}
  Hence, on $E_{0, \R^2}^1 \cap (\ker d_c)^\bot $, the operator
  $\delta_c d_c$ is conjugated by $d_{\R^2}$ to $P^*P = X^4 + Y^2$
  acting on the one dimensional bundle $\theta_Z$. Its spectrum can be
  described by Fourier analysis, in $\R^2$ here. The Fourier transform
  of the spectral space $E(P^*P \leq \lambda^2)$ are the $L^2$
  functions whose supports are contained in the set
  $$FS(\lambda) = \{(x, y) \in \R^2, x^4 + y^2 \leq \lambda^2\}. $$
  When $\lambda \rightarrow 0$, this strange collection of flying
  saucers are shrinking at speed $\lambda$ along $y$ and only
  $\sqrt\lambda$ along $x$. Therefore the rescaling $(X, Y)
  \rightarrow (\lambda^{-1/2} X, \lambda^{-1} Y )$ associated to the
  previous weight $N'$ appears naturally here. Observe anyway that the
  $FS(\lambda)$ have no usable self-similar limit in this rescaling
  (they converge to the union of two segments along $x$ and $y$). In
  fact one finds easily that the asymptotic area of $FS(\lambda)$ is
  $-\lambda^2 \ln \lambda$ when $\lambda \rightarrow 0$. This measure
  of $FS(\lambda)$ is also the density of $E(P^*P \leq \lambda^2)$,
  $E(\delta_c d_c \leq \lambda^2)$, and finally $E(\delta d \leq
  \lambda^2)$ as acting on $L^2(\R^2)$. We get finally by Laplace
  transform that the asymptotic $\R^2$-heat decay on $1$-forms is
  equivalent to $\dsp \frac{\ln t}{t^{1/2}}$ as $t \rightarrow
  +\infty$ (to be compared to the standard heat decay in $1/t$ on $\R^
  2$).

  Lastly, one can obtain some information on the anisotropic aspects
  of large scale heat diffusion here. For $\varepsilon > 0$, consider
  the orthogonal splitting $$L^2(\R^2) = L^2_\varepsilon (\R^2) \oplus
  (L^2_\varepsilon (\R^2))^\bot, $$
  with $L^2_\varepsilon (\R^2)$ the
  space of $L^2$ functions whose Fourier transforms are supported in
  $D_\varepsilon = \{(x,y), |y| \leq \varepsilon |x|\}$. We note that
  for each $\varepsilon$, $\mathrm{area} (FS (\lambda) \cap
  D_\varepsilon ) \asymp - \lambda\ln \lambda$ while $\mathrm{area}
  (FS (\lambda) \cap D_\varepsilon^\bot ) \asymp \lambda= o (-
  \lambda\ln \lambda) $, showing that most the spectral measure of
  $E(P^*P \leq \lambda^2)$ is supported on $D_\varepsilon$ when
  $\lambda \rightarrow 0$. Recall that the $d_c$ Laplacian on
  horizontal $1$-forms is conjugated to $P^*P$ acting on $\theta_Z$ by
  the map
  $$
  \delta_{\R^2} (f \theta_Z) = (Y.f) \theta_X -(X.f) \theta_Y .$$
  Since for functions $f \in L^2_\varepsilon (\R^2)$, one has
  $$\|Y.f\|_2 = \|y \hat f\|_2 \leq \varepsilon \|x \hat f\|_2 =
  \varepsilon \|X.f\|_2,$$
  we see that most the spectral measure of
  $E(\delta_c d_c \leq \lambda^2)$, and finally $E(\delta d \leq
  \lambda^2)$ is supported on forms $\alpha$ such that $\|\alpha
  (X)\|_2 \leq \varepsilon \|\alpha (Y)\|_2$, that is forms close to
  the $Y$-direction.  Topologically this direction corresponds to the
  largest Massey product we can find in the $1$-minimal model of $G$.
  Namely, one has $\theta_X$, $\theta_Y \in H^1(\mathfrak{g})$ with
  \begin{equation*}
    \left\{
      \begin{aligned}
        - \theta_X \wedge \theta_Y &  = d_0 \theta_Z\\
        - \theta_X \wedge \theta_Z & = d_0 \theta_T \\
        - \theta_X \wedge \theta_T & \in H^2 (\mathfrak {g}),
      \end{aligned}
     \right.
  \end{equation*}
  Analytically, this sequence is recognized in \eqref{eq:2:2:6} as the
  $X^3$ component of $d_c$ from $\theta_Y$ to $\theta_X \wedge
  \theta_T$, the link being done by taking the elliptic symbol of
  $d_c$ in the $X$ direction. This shows that some information on the
  multiplicative structure of $H^*(\mathfrak{g})$, not only its
  weight, is analytically encoded in $(E_0,d_c)$ even at the
  commutative level.  The effect here is that we can ``hear'' the
  exceptional direction $Y$.

\section{Towards a refined pinching}
\label{sec:5}

The pinching of $\alpha_1$ we gave relies on very few analysis,
general facts on the Hilbert-Schmidt and trace class operators, and
the ellipticity of the de Rham Laplacian. Anyway it is certainly not
sufficient to obtain the exact value of $\alpha_1$, $\beta_1$ on
groups of non homogeneous presentation. To improve (even partially)
the result we have to be more careful on the analytic properties of
$d_c$ as related to the convergence of the underlying spectral
sequence and filtered complex.

\subsection{C-C ellipticity and applications} 
\label{sec:5:1}

To that purpose we introduce some class of analytic regularity that
will fit to the vector valued operators on graded groups we are
dealing with. It relies on the notion of (maximal) hypoellipticity as
defined for scalar operators.

Let $G$ be a graded group, and $P$ a left invariant differential
operator on $C^{\infty}(G)$ of order $p$ with respect to natural
grading in $TG$.  $P$ is said maximally hypoelliptic if locally
$\|Pf\|_2 + \|f\|_2$ controls all derivatives of $f$ of weight $\leq
p$. One can show (see \cite{Ro-St,He-No2}) that this is equivalent to
the existence for $P$ of a parametrix $Q$ of order $- p$. That means
the kernel of $Q$ is an homogeneous distribution on $G$ of order
$p-N(G)$ near the diagonal, and such that $Q\circ P = \id +
\mathcal{S}$ for some smoothing operator $\mathcal{S}$. This amounts
also to the partial invertibility of $P$ at the convolution level in
$q * p = \delta_0 + s$ with $s\in \mathcal{S}(G)$ (where the
convolution product on $G$ is defined by $ (f*g) (x) = \int_G f(y^{-1}
x) g(y) dy$). This last criteria makes sense on more general operators
than differential ones. On a filtered group the basic example of a
second order maximally hypoelliptic operator is given by the
horizontal sub-Laplacian defined on functions by
\begin{equation}\label{eq:5:1}
  \Delta_H f = - \sum_{i=1}^k X_i^2 f, 
\end{equation}
where $D=\mathrm{span} (X_i, 1\leq i\leq k)= \mathfrak{g}_1$ generates
$\mathfrak{g}$. Starting from this positive self-adjoint operator, one
can define $|\nabla| = \Delta_H^{1/2}$ which can be shown to be an
order one hypoelliptic pseudo-differential operator on $G$ (see
\cite{Fo} section 3) invertible in the previous sense. On a general
graded group, not necessarily filtered, such an order one scalar
hypoelliptic operator can also be constructed. One can either replace
$\Delta_H$ with the operator $P$ of \cite{He-No} corollaire 0.2, or
use the one of \cite{Ch-Ge} section 6. We denote by $|\nabla|^{-1}$
its parametrix, in both left and right senses by self-adjointness of
$|\nabla |$.

We are now ready to define the class of regularity we need. We extend
$|\nabla|$ as acting diagonally on the full exterior algebra $\Omega^*
G$. We call $N$ the induced weight function there.
\begin{defn}
  An operator $P$ acting on (part) of $\Omega^* G$ is called C-C
  elliptic if $P^{\nabla} = |\nabla|^{-N} P |\nabla|^{N}$ is a
  maximally hypoelliptic operator in the previous sense.
\end{defn}
In the scalar case this is just hypoellipticity. Thus $P=\Delta_H$ is
such an operator, while the full de Rham Laplacian on functions is not
(except in $\R^n$).  Indeed if $G$ is a $r$-step group, $\Delta$ is a
$2r$-order operator in the graded sense, as it contains $X_r^2$ terms
with $X_r \in \mathfrak{g}_r$, but $\Delta$ do not control
$2r$-derivatives along $\mathfrak{g}_1$ for $r > 1$. One interesting
feature of this class is that even if de Rham Laplacian is not in it,
the de Rham complex itself is, as $(E_0, d_c)$.
\begin{thm}\cite{Ru}
  \label{thm:5:2}
  Let $G$ be a graded nilpotent Lie group. Then the de Rham complex
  and $(E_0, d_c)$ are C-C elliptic.
\end{thm}
That means the Laplacians associated to $d^\nabla$ and $d_c^\nabla$
\begin{equation}
  \Box_d = d^\nabla
  (d^\nabla)^* + (d^\nabla)^* d^\nabla\ \mathrm{and}\ \Box_{d_c} =
  d_c^\nabla (d_c^\nabla)^* + (d_c^\nabla)^* d_c^\nabla
\end{equation}
(\emph{not} $(\Delta)^\nabla$ and $(\Delta_c)^\nabla$) are maximally
hypoelliptic. Observe that they are all $0$ order operators !  Indeed
$d$ and $d_c$ are homogeneous operators on $\Omega^*G$ in the sense
the components of $d$ and $d_c$ that increase the weight of forms by
$k$ are precisely order $k$ operators (the whole construction of
section \ref{sec:2} is homogeneous in that sense).
\begin{proof}
  Thanks to analysts' work a short proof of thm \ref{thm:5:2} is
  possible.
  
  It relies on an hypoellipticity criteria (called Rockland's
  condition) first proved for differential operators on general graded
  groups by Helffer and Nourrigat \cite{He-No}, and extended to the
  pseudo-differential setting by Christ and all in \cite{Ch-Ge}. Given
  a (pseudo-differential) operator $P$ of order $k$ on $G$, call $P_k$
  its homogeneous part of order $k$. Then the maximal hypoellipticity
  of $P$ is implied by the injectivity, for all non trivial
  irreducible unitary representations $\pi$ of $G$, of $\pi (P_k)$ on
  $\mathcal{S}_\pi \subset H_\pi$ the space of $C^\infty$ vectors of
  the representation. (If $G =\R^n$ this is the well known ellipticity
  criteria on the injectivity of the principal symbol of $P$ in all
  directions.  In general we recall (see \cite{Co-Gr} or
  \cite{He-No2}) that Kirillov has shown that these representations
  are parameterized by non-zero covectors $\xi \in \mathfrak{g}^*$
  modulo the co-adjoint action of $G$ on $\mathfrak{g}^*$.)
  
  We will need to know that for $0 \not= \pi \in \widehat{G}$ one can
  find some $X\in \mathfrak{g}$ such that $\pi (X)$ is the scalar
  (anti-self adjoint) operator $i \id$.  Indeed, by irreducibility of
  $\pi$ the center $Z(G)$ of $G$ acts by scalars on $H_\pi$ (given by
  a character $Z(G)\rightarrow U(1)$). If it is trivial, consider the
  quotiented action on $Z(G/Z(G))$ and so on until, by nilpotency of
  $G$, some non trivial scalar action is found.
  
  We come to the proof. Let $\pi$ and $X$ as above.  that the $0$
  order part of $\pi(\Box_d)$ and $\pi (\Box_{d_c})$ are injective on
  the $C^\infty$ vectors $\mathcal{S}_\pi \subset H_\pi$.  These
  spaces really identify with the standard Schwarz space
  $\mathcal{S}(\R^n) \subset L^2(\R^n)$ with $\R^n \simeq G/H$ when
  $\pi$ is realized as an induced representation from $G/H$. Invariant
  differential operators on $G$ become differential operators with
  polynomial coefficients on $\R^n$, while by definition $|\nabla|$ is
  such that $\pi (|\nabla|)$ is invertible and preserves
  $\mathcal{S}_\pi$ (in the filtered part we have just $\pi (|\nabla|)
  = \pi (\Delta_H)^{1/2}$).  Hence the injectivity of $\pi(\Box_d)$ on
  $\mathcal{S}_\pi$ is equivalent to the injectivity of the system
  $\pi (d^\nabla) + \pi (d^\nabla)^*$ there.  This is certainly
  implied by
  \begin{displaymath}
    \ker \pi (d^\nabla) \cap \mathcal{S}_\pi = \pi (d^\nabla)
    (\mathcal{S}_\pi).
  \end{displaymath} 
  This is equivalent by dropping the conjugation by $|\nabla|^N$,
  which preserves $\mathcal{S}_\pi$, to the vanishing of the smooth
  cohomology of the complex $\pi (d)$
  \begin{equation}
    \label{eq:5:3}
    \ker \pi (d) \cap \mathcal{S}_\pi = \pi (d) (\mathcal{S}_\pi).
  \end{equation}
  Consider $\mathcal{L}_X$ the Lie derivative along $X$, with $X$ as
  above. As a differential operator on $\Omega^*G$ it can decomposed
  in $\mathcal{L}_X = X \id + ad(X)$ where $$ad(X)\alpha = - \sum
  \alpha(\cdot, [X, \cdot], \cdot) = (\mathcal{L}_X)_0 = i_X d_0 + d_0
  i_X$$
  is the algebraic part of $\mathcal{L}_X$ when expressed in a
  left invariant base.  Since $\pi (X) = i \id$, one has $\pi
  (\mathcal{L}_X) = i\id + ad(X)$ on $H_\pi$ with $ad(X)$ nilpotent.
  It is therefore invertible of inverse
  $$
  P_X = \pi (\mathcal{L}_X)^{-1} = \sum_{k \geq 0} i^{k-1}
  ad^k(X).$$
  From Cartan's formula $\mathcal{L}_X = i_X d +di_X$ we
  get $[\mathcal{L}_X , d] = [\mathcal{L}_X, i_X] =0$, so that $[P_X ,
  \pi (d)] = [P_X, i_X] = 0$, and finally
  \begin{displaymath}
    \id = (P_X i_X) \pi (d) + \pi (d) (P_X i_X),    
  \end{displaymath}
  which gives \eqref{eq:5:3}. The proof for $(E_0, d_c)$ is a formal
  consequence of this one. We recall that by theorem \ref{thm:2:2:1},
  $d_c= \Pi_{E_0} \Pi_E d \Pi_E \Pi_{E_0}$, with $\Pi_{E_0} \Pi_E
  \Pi_{E_0} = \id $ on $E_0$. Composing \eqref{eq:5:4} with $\pi
  (\Pi_{E_0} \Pi_E)$ on the left and $\pi (\Pi_E \Pi_{E_0})$ on the
  right we get on $E_0$
  \begin{displaymath}
    \id = Q_X \pi (d_c) +  \pi (d_c) Q_X,
  \end{displaymath}
  with $Q(X) = \pi (\Pi_{E_0} \Pi_E) P_X i_X \pi (\Pi_E \Pi_{E_0})$.
  As before this implies the injectivity of $\pi(\Box_{d_c})$.
\end{proof}
\begin{rem}
  The role of the conjugation by $|\nabla|^N$ in this proof appears
  very formal. Anyway it is crucial to put everything in a Sobolev
  scale and work with an homogeneous operator.
\end{rem}

The twisted $d_c^\nabla$ has a rather strange harmonic theory, due to
the fact that $(d_c^\nabla)^* = |\nabla|^N \delta_c |\nabla|^{-N}$ is
not $\delta_c^\nabla = |\nabla|^{-N} \delta_c |\nabla|^N$. On the
other hand we observed in remark \ref{rem:3:7} that the $d_c$-harmonic
theory is not so attractive, because not scale invariant, except in
the case we have $E_0^k = H^k (\mathfrak{g})$ of pure weight $N_k$.
This is precisely the assumption we did in the previous pinching
results. In that case the C-C regularity $d_c$ is actually useful for
the $d_c$ harmonic theory itself. We remark that C-C ellipticity is
preserved by any conjugation by an invertible fixed $|\nabla|^{N_0}$.
So one can use any shifted weight function $N-N_0$ instead of $N$. In
the case $E_0^k$ is of pure weight $N_k$ the interesting one is $N' =
N-N_k$, since then on $E_0^k$
\begin{equation}\label{eq:5:4}
  d_c^\nabla + (d_c^\nabla)^* = |\nabla|^{-N'} d_c + |\nabla|^{N'}
  \delta_c = P (d_c + \delta_c),
\end{equation}
with $P$ of \emph{strictly} negative order $\leq -1$, because $d_c$
increases the weight by $1$ at least.

We now show some applications of these techniques. Let again $G$ be a
graded group $G$ and some $k$ with $E_0^k = H^k (\mathfrak{g})$ of
pure weight $N_k$. The C-C ellipticity of $(E_0, d_c)$ allows to
define some number $r_k$ describing the analytic order of convergence
of the spectral sequence associated to $d$ on $k$-forms.  Starting
form the splitting \eqref{eq:3:3:1} of $d_c$ in $d_c^{\delta
  N_k^\mathrm{min}} + \cdots + d_c^{\delta N_k^\mathrm{max}}$, we
consider for $r \in [N_k^\mathrm{min}, N_k^\mathrm{max}]$ its cut-off
\begin{displaymath}
  d_c^{[r]} = d_c^{\delta N_k^\mathrm{min}} + \cdots + d_c^r.
\end{displaymath}
By theorem \ref{thm:5:2}, we know that the complex
\begin{displaymath}
  \xymatrix{E_0^{k-1} \ar[r]^{d_c} & E_0^k \ar[r]^{d_c^{[r]}} & E_0^{k+1}
    }
\end{displaymath}
is C-C elliptic on $E_0^k$ if $r = \delta N_k^{\mathrm{max}}$.
\begin{defn}\label{def:5:3}
  Define $r_k$ as the smallest integer $r$ such that the previous
  complex is C-C elliptic on $E_0^k$.
\end{defn}
In view of the proof of theorem \ref{thm:5:2} and \eqref{eq:5:3} this
$r_k$ is the smallest number $r$ such that $\ker \pi (d_c^{[r]}) \cap
\mathcal{S}_\pi = \pi (d_c) (\mathcal{S}_\pi)$ for all non-trivial
$\pi \in \widehat{G}$. From a purely algebraic viewpoint, this $r_k$
is actually the order of convergence, at the level of representation,
of the underlying spectral sequence associated to the natural filtered
complex here (see section \ref{sec:2}). In other words, $r_k$ is also
the smallest integer such that for all $\pi$, any smooth $k$-form
$\alpha \in \mathcal{S}_\pi$ of weight $p$ which satisfies $\pi(d)
\alpha$ of weight $\geq p+ r$ is necessarily $\pi(d)$-closed (and
finally exact) up to a form of weight $p+1$.
\begin{thm}
  \label{thm:5:4}\cite{Ru} Let $G$, $E_0^k$ as in theorem
  \ref{thm:3:3:1}. Then for $r_k$ as above
  \begin{displaymath}
    \delta
    N_k^\mathrm{min} \leq \beta_k \leq r_k \leq \delta N_k^\mathrm{max}.
  \end{displaymath}
\end{thm}
We will give examples of groups with $r_k < \delta N_k^\mathrm{max}$
in the next section.
\begin{proof}
  We have to show that $\beta_k \leq r_k$. We use the same
  near-cohomology approach as in the proof of theorem \ref{thm:3:3:1}.
  Firstly since for $\alpha \in E_0^k$
  $$
  \|d_c^{[r_k]} \alpha\|^2 = \|d_c^{\delta N_k^\mathrm{min}}
  \alpha\|^2 + \cdots + \|d_c^{r_k} \alpha\|^2 \leq \|d_c \alpha\|^2,
  $$
  we certainly have that the near-cohomology cones thickness
  functions of $d_c$ and $d_c^{[r_k]}$ satisfy
  $$F_{d_c} (\lambda) \leq F_{d_c^{[r_k]}} (\lambda).$$
  Hence the same
  homogeneity argument as in the proof of thm \ref{thm:3:3:1} gives
  the result if we show that $\dsp F_{d_c^{[r_k]}} (\lambda)$ is
  finite. This is a consequence of the C-C ellipticity assumption.
  Namely the same proof as theorem \ref{thm:5:2} together with
  \eqref{eq:5:4} shows that $ |\nabla|^{-N'} d_c^{[r_k]} +
  |\nabla|^{N'} \delta_c = P (d_c^{[r_k]} + \delta_c )$ is an order
  $0$ maximally hypoelliptic pseudo-differential operator. That means
  (see \cite{Ch-Ge}) it is bounded and invertible in $L^2$ up to a
  smoothing term
  \begin{displaymath}
    QP (d_c^{[r_k]} +  \delta_c ) = \id + \mathcal{S},
  \end{displaymath}
  for some bounded (order $0$) $Q$ and smoothing $\mathcal{S}$. We
  have already observed that $P$ is of order $\leq -1$, hence $B=
  |\nabla| QP$ is bounded and satisfies
  \begin{equation}
    B (d_c^{[r_k]} +  \delta_c ) = |\nabla| + \mathcal{S}.
  \end{equation}
  Therefore there exists some $C$ such that
  \begin{displaymath}
    \||\nabla| \alpha\| \leq C (\|d_c^{[r_k]} \alpha\| + \|\delta_c
    \alpha\| + \|\alpha\|),  
  \end{displaymath}
  giving on $\ker \delta_c$ the control of $F_{d_c^{[r_k]}} (\lambda)$
  by $F_{|\nabla|} (C \lambda + C)$. But this last function is finite
  because $|\nabla|$ is a first order maximally hypoelliptic operator
  with smoothing spectral projectors (for any $p$, $|\nabla|^{-k}$ has
  a $C^p$ kernel for $k$ large enough by the Sobolev embedding
  theorems).
\end{proof}

In the opposite direction to the previous result, one can sometimes
improve the bound from below of $\beta_k$, if $d_c^{[r]}$ is
sufficiently degenerated.
\begin{prop}
  \label{prop:5:5} Let $G$, $E_0^k$ as above. Suppose that for some
  $r$, there exists non-closed forms in $D(d_c) \cap \ker d_c^{[r]}$,
  then $\beta_k \geq r + 1$.
\end{prop}
\begin{proof}
  Consider the closed quadratic form $q(\alpha) = \|d_c \alpha\|^2$ on
  $(\ker d_c)^\bot$. Its restriction to $D(d_c) \cap \ker d_c^{[r]}$
  is still closed, positive, and therefore associated to a
  self-adjoint $\Delta_r$ (see \cite{Re-Si} thm VIII.15). By
  hypothesis $\Delta_r \not= 0$. Then for some $\lambda$, $L = E(0<
  \Delta_r \leq \lambda^2)$ is a non zero $G$-invariant linear space
  in the near-cohomology cone $C_\lambda(d_c)$. Since $d_c$ is
  homogeneous of order $\geq r+1$ on $D(d_c) \cap \ker d_c^{[r]}$ we
  have, as in the proof of theorem \ref{thm:3:3:1}, that
  $h_\varepsilon^* (L) \subset C_{\lambda \varepsilon^{r+1}}(d_c)$.
  This gives the bound from below for $\beta_k$.
\end{proof}
\begin{rem}
  We have already used this for the triangular group $N_4$ in section
  \ref{sec:4:2}.
\end{rem}

We can interpret the assumption in proposition \ref{prop:5:5} in terms
of representation theory. Namely by the Fourier-Plancherel isomorphism
between $G$ and $\widehat{G}$, the existence of some non-zero $L^2$
form $\alpha \in D(d_c) \cap \ker d_c^{[r]}$ implies that $\pi
(\alpha)$ is a field of Hilbert-Schmidt operators in $\ker \pi
(d_c^{[r]})$ with $\pi(\alpha)\not= 0$ on a set of representations of
strictly positive Plancherel measure. Then the cohomology $\ker \pi
(d_c^{[r]}) / \im \pi (d_c)$ is non-vanishing for ``a lot'' of
representations of maximal dimensions (called generic). This is
opposite (but not complementary) to the case of definition
\ref{def:5:3} and theorem \ref{thm:5:4}, where we required this
cohomology to vanish on all representations.

We close this section with an application the $d_c$ Hodge-de Rham
decomposition on compact $E_0$ regular C-C manifolds $M$. We have seen
in proposition \ref{prop:3:6} that the we have a closed splitting in
$L^2(\Omega^*M)$
\begin{displaymath}
  E_0 = \mathcal{H}_c \oplus \im d_c \oplus \im \delta_c,
\end{displaymath}
where $\mathcal{H}_c = \ker d_c \cap \ker \delta_c$.  Again we
restrict to the case of $E_0^k=H^k(\mathfrak{g}_{x_0})$ is of pure
weight $N_k$ (see also remark \ref{rem:3:7}). We assume also that C-C
structure is given by a bracket generating distribution $D$ (the
filtered case).
\begin{prop}\label{prop:5:7}
  Let $M$, $E_0^k$ as above. Then $\mathcal{H}_c = \ker d_c \cap \ker
  \delta_c$ consists of smooth forms.
\end{prop}
To see this we define a differential operator of the form $\Box_c = P
(d_c + \delta_c)$, with $P$ chosen such that $\Box_c$ becomes
homogeneous. Fix any horizontal connection $\nabla_D$ on $\Omega^*M$.
let $K$ be the maximal order of the components of $d_c + \delta_c$ on
$E_0^k$. Then we consider
\begin{displaymath}
  \Box_c = \sum_p \nabla_D^{K-p} (d_c^p + \delta_c^p),
\end{displaymath}
where $d_c^p$ (resp. $\delta_c^p$) is the components of $d_c$ that
increases (resp. decreases) the C-C weight by $p$. This $\Box_d$ is a
differential operator of order $K$. We show that
\begin{lemma}\label{lemma:5:9}
  $\Box_c$ is maximally hypoelliptic on $M$ (in the bracket generating
  case).
\end{lemma}
We recall this means that $\Box_c \alpha$ controls the full the
$K$-horizontal jet of $\alpha$ (in $L^2$ norms). One consequence of
this is that weak solutions to $\Box_c \alpha =0$ are smooth, giving
proposition \ref{prop:5:7}.
\begin{proof}
  Some characteristic features of maximal hypoellipticity of
  differential operators are its stability under perturbations of
  lower order and the fact, in the filtered case, it can be checked on
  the model group $G_{x_0}$ with the freezed operator $\Box_{c,x_0}$
  (see \cite{He-No2}). We have already observed that $d_c$ at $x_0$
  may be viewed as a perturbation of $d_c$ on $G_{x_0}$ (by Cartan's
  formula $d\alpha (X_i) = \sum X_i \alpha (\cdot) - \alpha (\cdot
  ,[X_i,\cdot ],\cdot )$ and the fact that $[\ , \ ]_0$ is the lower
  part of $[\ ,\ ]$).  Hence we just have to check the result on the
  tangent group. There $\Box_c$ looks like an algebraic version of
  $d_c^\nabla + (d_c^\nabla)^*$ in \eqref{eq:5:4} with the shift
  weight $N' = N-N_k$, except the scalar pseudo-differential
  $|\nabla|$ is replaced by the full first order horizontal jet
  $\nabla_D$. The same proof as in theorem \ref{thm:5:2} actually
  gives the injectivity of $\pi (\Box_c)$ since on $\mathcal{S}_\pi$,
  $\ker \pi (\Box_c) =\ker \pi (d_c) \cap \ker \pi (\delta_c)$ and we
  have already seen that $\ker \pi (d_c) \cap \mathcal{S}_\pi = \pi
  (d_c) (\mathcal{S}_\pi)$.
\end{proof}
\begin{rem}
  The last part of this proof applies on any graded group, not only
  the filtered ones. It can therefore be taken as an alternative
  approach to C-C ellipticity using only differential operators, in
  this case of $E_0^k$ of homogeneous weight.
\end{rem}

\subsection{Last examples}
\label{sec:5:2}

We start with a series of groups illustrating the ``analytic''
pinching of $\beta_1$ obtained in theorem \ref{thm:5:4}.  We show that
some groups may have (arbitrarily) high order relations and still
$r_1$ and thus $\beta_1 = 1$.

\subsubsection{First inaudible relations }

Consider some quadratically presented group $G$ of rank $r>3$. The
examples of Carlson-Toledo and Chen have been described in section
\ref{sec:4:1}. Let $I$ be any ideal in $\mathfrak{g}$ generated by
elements of weight $\geq 3$, $N = \exp I$ and consider $H= G/ N$. The
relations of $H$ as referred to the free Lie group, are generated by
the quadratic relations of $G$ and the generators of weight $\geq 3$
of $N$. It may therefore have high order relations, but anyway they
are ``inaudible'' in the heat decay.
\begin{prop}
  For such groups $H$, $r_1= 1= \beta_1 (H)$, and finally $\alpha_1
  (H) = N(H)$.  In particular these values are the same as on
  quadratic groups.
  \end{prop}
  \begin{proof}
    We use theorem \ref{thm:5:4}. Let $\pi_H$ be a non-trivial
    irreducible unitary representation of $H$. Suppose $\alpha$ is a
    form in $\mathfrak{h}_1^* \cap \mathcal{S}_\pi$ such that $\pi_H
    (d_c^{[1]}) \alpha =0$.  This means that $\alpha$ has an extension
    $\widetilde {\alpha}$ (given by $\pi_H (\Pi_E) \alpha$) such that
    $\pi_H (d) \widetilde {\alpha}$ is of weight $\geq 3$.  Consider
    $P : G \rightarrow H$. Then $P^* (\pi_H)= \pi_G$ is an irreducible
    unitary representation of $G$ and $\pi_G (d) (P^* \widetilde
    {\alpha}) = P^* (\pi_H (d) \widetilde{\alpha})$ is still of weight
    $\geq 3$. But since $G$ is quadratically presented this implies
    that $\pi_G (d_c^G) (P^*\widetilde {\alpha}) = 0$, and finally
    $P^* \widetilde {\alpha}$ can be written $\pi_G (d_c^G) f$ for
    some $f \in \mathcal{S}_\pi$.  Restricting this to
    $\mathfrak{g}_1$, gives $\alpha = \pi_F (d_c) f$ and the result.
  \end{proof}
  
  These examples show that the spectral sequence may converge in the
  $L^2$ sense quicker than in the algebraic one. More precisely $H^2
  (\mathfrak{h})$ always contain forms of weight $\geq 3$ here, dual
  to the relations we added by proposition \ref{prop:2:9}.  View as
  left invariant forms in $H$ they are closed, therefore locally
  exact. Thus we obtain one forms $\alpha$ such that $d\alpha \in H^2
  (\mathfrak{h})$ are non-zero forms of weight $\geq 3$. In other
  words $\alpha_D = \Pi_{E_0}\alpha$, the restriction of $\alpha$ to
  $D= \mathfrak{h}_1$ satisfies $d_c \alpha_D$ of weight $\geq 3$.
  This shows, as claimed, that the convergence rank of the spectral
  sequence as a local tool is actually greater than in $L^2$.
  
  We can explain (partially) more geometrically why these high order
  relations do not have an $L^2$ trace. We have to study
\begin{equation}\label{eq:E2}
  E_2 = \{ \alpha \in E_0^1(H) = \Omega^1 D \mid d_c
  \alpha \ \mathrm{is\ of \ weight\ }\geq 3\}.
\end{equation}
\begin{prop}\label{prop:5:2:1}
  Up to closed forms, $E_2$ contains only forms whose components are
  polynomial functions (of bounded degree) in the coordinates of $G$.
\end{prop}
Thus certainly $E_2$ do not contain any non-closed $L^2$ form. We
can't make an $L^2$ ``wave packet'' of them.
\begin{proof}
  We first show than the components of $d_c \alpha$ are polynomial
  functions. Then lifting $\alpha$ to $\Omega^1 H$ (with $\pi_E$) and
  integrating $d \widetilde{\alpha}$ with Poincar\'e's lemma along
  polynomial vector fields gives the result.
  
  Let $\Pi: G \rightarrow H$ and $\alpha \in E_2$. Again $\Pi^*\alpha$
  is in $E_2(G)$ which are $d_c^G$-closed because $G$ is quadratically
  presented. That means $\Pi^*\alpha$ has a true closed extension
  $\beta$ on $G$. Recall that by section \ref{sec:2:4}, $d_c \alpha$
  interprets as the components of $\beta$ along the generators of $I$
  the ideal of relations of $H$ with respect to $G$. So we have to
  show these components are polynomial. This is consequence of
  \eqref{eq:2:4:3}. Namely, suppose $G$ is of rank $r$. Applying
  \eqref{eq:2:4:3} to $Y \in I$ of weight $r$, gives $X.\beta(Y)=0$
  because $[X,Y]=0$ for any $X$. Therefore $\beta(Y)$ is constant.
  Applying now \eqref{eq:2:4:3} to $Y\in I$ of weight $r-1$ gives then
  that $X.\beta (Y)$ are constants functions, since $[X,Y]$ is in $Y$
  and of weight $r$. And so on, all components of $\beta$ along $I$
  are polynomials.
\end{proof}

\subsubsection{Around $H$-groups}

The previous groups give us first example of ``inaudible'' relations.
We give another illustration of the fact that ``audible'' relations
have to be sufficiently flexible. We consider again $H$-type groups
(see \ref{sec:4:1}). A complete study may be found in \cite{Ka} or
\cite{Co-Do}. Let $D=\R^n$, and $\R^k$ be endowed with scalar
products. A Clifford module structure on $D$ is a linear map
\begin{displaymath}
  J: \R^k \rightarrow GL(n)
\end{displaymath}
such that $J(\theta)^2 = - \|\theta\|^2 \id$ and $J(\theta)^* =
-J(\theta)$.  Consider then
\begin{align*}
  L : \R^k & \rightarrow \Lambda^2 D^* \\
  \theta & \rightarrow g_D( J(\theta) \cdot, \cdot).
\end{align*}
The $H$-type group associated to $J$ is the two step group $G$
determined by the extension of $D$ by $\im L$. That means
$\mathfrak{g}= D \oplus T$ with the Lie bracket structure $D\times D
\rightarrow T=\R^{k*}$ given by $\theta([X, Y]) = -L(\theta) (X,Y) = -
\langle J(\theta ) X, Y\rangle$. When $k=1, 3, 7$ such groups appears
as the maximal nilpotent (Iwasawa) subgroups of the rank 1 semi-simple
groups $SU(n,1)$, $Sp(n,1)$ and $F_4^{-20}$ (acting on the Cayley
plane).
\begin{prop}
  Let $G$ be a $H$-type group, then either $\beta_1 (G) =1$ or
  $\beta_1 (G) = 2$ and $G$ is the 3 dimensional Heisenberg group, or
  the 7 dimensional quaternionic group (associated to an elliptic $D^4
  \subset \R^7$ see section \ref{sec:4:1}).
\end{prop}
\begin{proof}
  Observe that $D$ is necessarily of even dimension $n=2p$, and that
  the $p=1$ case corresponds to to the Heisenberg group $\mathbf{H}^3$
  we have already met ($\beta_1= 2$ by cubic presentation). So we
  restrict now to $p\geq 2$. We first show how the use of theorem
  \ref{thm:5:4} reduces things to a problem on the Heisenberg group
  $\mathbf{H}^{2p+1}$. We study $r_1$. By the orbit method non-trivial
  irreducible representation of $G$ are of two types, up to the
  coadjoint action.
    
\item - The first ones are trivial on $T$ and come from some character
  $\pi : D \rightarrow U(1)$ with $\pi (V) = e^{i \langle \xi
    ,V\rangle}$ and non zero $\xi \in D$. We show that $\ker (\pi
  (d_c^{[1]})) \subset \im \pi (d_c)$ on $\Lambda^1 D^*$.  Namely $\pi
  (d_c^{[1]}) \alpha =0$ means that $i \xi \wedge \alpha \in \im L$
  (because $d_c^{[1]}$ is the quotient of the horizontal differential
  by $\im L$).  But $\im L$ contains only non-degenerate forms of top
  rank $2p \geq 4$ here. Therefore $\xi \wedge \alpha =0$ and $\alpha
  = iC \xi = \pi (d_c) C$.
    
\item - The other (generic) ones are induced by a vertical $\theta$,
  and factors trough a representation $\pi$ of the Heisenberg group
  $\mathbf{H}^{2p +1}$ with the contact form $\theta$ and $d \theta =
  g_D( J(\theta) \cdot, \cdot)$. We orthogonally split $\im L= \R
  d\theta \oplus L'$, and observe that by the Clifford structure $L'
  \subset \Lambda^{2,0}D^* \oplus \Lambda^{0,2}D^*$ with respect to
  $J(\theta)$ (while $d\theta \in \Lambda^{1,1} D^*$). Again the
  equation $\pi (d_c^{[1]}) \alpha = 0$ means $\pi(d_D) \alpha \in \im
  L$, and our system can finally be interpreted on $\mathbf{H}^{2p
    +1}$, as
    \begin{equation}\label{eq:5:2:1}
      \pi(d_c^\mathbf{H}) \alpha \in L' \subset \Lambda^{(2,0) +
      (0,2)}D^* \quad \mathrm{and}\quad \pi(\delta_c^\mathbf{H})
      \alpha= 0,
    \end{equation}
    with $\alpha \in \mathcal{S}_\pi$ a smooth vector of the
    representation. We arrived at the turning.
    \begin{lemma}\label{lemma:5:2:1}
      Equation \eqref{eq:5:2:1} has a non trivial solution iff there
      exists a non zero $\alpha \in D^*$ such that $$\alpha^{1,0}
      \wedge \Lambda^{1,0} D^* \subset L^{'(2,0)}.$$
      In real notations
      $\alpha \wedge \beta - J \alpha \wedge J\beta \in L'$ for all
      $\beta \in D^*$.
    \end{lemma}
    Note that in this case the forms $\alpha \wedge \beta - J \alpha
    \wedge J\beta \in L' \subset \im L$ are of rank 4. Since here $\im
    L$ contains only top rank forms one has necessarily $D$ of
    dimension $4$.  Moreover counting dimensions shows that
    $\mathrm{dim} D - 2 \leq \mathrm{dim} (\im L) -1 = k-1$, ie $k
    \geq 3$, and finally $k=3$ since this is the maximal number of
    complex structures on $D=\R^4$. Then certainly $G$ has to be the
    seven dimensional quaternionic group $Q_7$. If not
    \eqref{eq:5:2:1} has no solution and theorem \ref{thm:5:4} gives
    that $r_1=1$ and thus $\beta_1= 1$.
    
    Reversely $Q_7$ satisfies the condition of lemma
    \ref{lemma:5:2:1}, because $L^{'(2,0)} =\Lambda^{(2,0)} D^* = \C
    dZ_1 \wedge dZ_2$. It remains to see that $\beta_1 (Q_7) =2$. By
    proposition \ref{prop:5:5} it is sufficient to find some non
    $d_c$-closed $\alpha \in \Omega^1 D \cap D(d_c)$ such that
    $d_c^{[1]} \alpha =0$. This is achieved by making a wave packet of
    solutions $\alpha_\pi$ of \eqref{eq:5:2:1}. Namely fix some
    compact set $K$ of positive Plancherel measure in the generic
    representations of $Q_7$. Let $P_{\alpha_\pi}$ be the orthogonal
    projection on $\alpha_\pi \subset \mathcal{S}_\pi$. Fourier
    inversion formula \eqref{eq:fourier} leads to consider for $g\in
    Q_7$
    \begin{displaymath}
      \alpha (g) = \int_{K \subset \widehat G} \tr (\pi (g)
      P_{\alpha_\pi} )  d\mu (\pi) = \int_K \langle \pi (g) \alpha_\pi
      , \alpha_\pi 
      \rangle_{H_\pi} d\mu (\pi). 
    \end{displaymath}
    $\alpha_\pi$ being smooth vectors, $g\rightarrow \alpha (g)$ is a
    smooth form. Moreover
    $$\|\alpha_\pi\|^2 = \int_K \|P_{\alpha_\pi}\|_{HS}^2 d\mu (\pi) =
    \mu (K) \in ]0,+\infty[,$$
    and for any derivative $\|X^I
    \alpha\|^2 = \int_K \|\pi (X^I) \alpha_\pi\|^2 d\mu (\pi) <
    \infty$. Of course $d_c^{[1]}\alpha = \int_K \langle \pi
    (d_c^{[1]}) \alpha_\pi, \alpha_\pi \rangle d\mu(\pi) =0$ while
    $d_c \alpha \not= 0$ because $\pi(d_c) \alpha_\pi \not= 0$ by C-C
    ellipticity of $d_c$.
  \end{proof}
  We are left with the proof of lemma \ref{lemma:5:2:1}
\begin{proof}
  We drop the $\pi$ and $\mathbf{H}$ to lighten notations and work on
  $\mathbf{H}^{2p+1}$ (for $p\geq 2$) with the contact complex $d_c$.
  A know fact is that the Laplacian $\Delta_c = (p-k) \delta_c d_c+
  (p-k+1) d_c \delta_c$ preserves the bigrading of $\Lambda^k D^*$
  (for $k< p$).  Applying this to $\Delta_c \alpha$ gives here
  \begin{align*}
    (\Delta_c \alpha,\alpha) & = (p-1) \|d_c^{2,0} \alpha^{1,0} \|^2 +
    (p-1) \|d_c^{0,2} \alpha^{0,1} \|^2\ \mathrm{by}\ \eqref{eq:5:2:1} \\
    & = (\Delta_c \alpha^{1,0}, \alpha^{1,0}) + (\Delta_c
    \alpha^{0,1}, \alpha^{0,1}) \\
    & = p \|\delta_c \alpha^{1,0}\|^2 + (p-1) \|d_c^{2,0} \alpha^{1,0}
    \|^2 +
    (p-1) \|d_c^{1,1} \alpha^{1,0} \|^2 \\
    & \quad + p \|\delta_c \alpha^{0,1}\|^2 + (p-1) \|d_c^{0,2}
    \alpha^{0,1} \|^2 + (p-1) \|d_c^{1,1} \alpha^{0,1} \|^2,
  \end{align*}
  leading to the vanishing of $\delta_c \alpha^{1,0}$, $d_c^{1,1}
  \alpha^{1,0}$, $\delta_c \alpha^{0,1}$ and $d_c^{1,1} \alpha^{0,1}$.
  By \cite{Ru2}, this is equivalent to the holomorphy of each
  component of $\alpha^{1,0}$, meaning $\overline{Z_i}.\alpha (Z_j)
  =0$ for $Z_{i,j} \in D^{1,0}$ (resp. anti-holomorphy of
  $\alpha^{0,1}$). At the representation level holomorphic functions
  are generated by $\dsp f=e^{-\sum_{i=1}^p x_i^2/2} \in
  \mathcal{S}(\R^p) = \mathcal{S}_\pi$, the vacuum state of the
  harmonic oscillator. Hence there exists a fixed $\beta^{1,0}$ such
  that $\alpha^{1,0}(x) = f(x) \beta^{1,0}$. Differentiating gives
  \begin{displaymath}
    \pi (d_c) \alpha^{1,0} = \pi (d_c^{2,0}) \alpha^{1,0} =  \sum
    (\pi (Z_i) f) \theta_i \wedge \beta^{1,0} \in L^{'(2,0)},
  \end{displaymath}
  but the functions $\pi (Z_i) f= (\frac{\partial}{\partial x_i} -
  x_i) f = - 2x_i f$ are independent giving that each $\theta_i \wedge
  \beta^{1,0}$ belongs to $L^{'(2,0)}$.
\end{proof}

We observe that this result gives rise to a group in the $\beta_1= 1$
class not quadratically presented, and yet of a different type than in
the previous section. Consider the 6 dimensional $G_6$, quotient of
$Q_7$ by a direction of its center, say $T_3$ associated to $\theta_3,
J_3$.  This $G_6$ is the tangent group to a generic $D^4 \subset
\R^6$.  It is a two dimensional extension of $D=\R^4$ with 2
orthogonal complex structures $J_1$ and $J_2$.

Consider again $E_2$ as in \eqref{eq:E2}. This $E_2$ has a lot of
non-closed sections. For instance in view of the proof of lemma
\ref{lemma:5:2:1}, it contains any form $\alpha \in \Lambda^1 D^*$
whose $(1,0)$ components with respect to $J_3$ are $J_3$ holomorphic
functions on $D$ invariant along $T_1, T_2$. Here $E_2$ contains many
other forms than polynomials of bounded degree like in proposition
\ref{prop:5:2:1}.  Thus the cubic relations of $G_6$ (dual to the
forms $\theta_1\wedge J_1\alpha - \theta_2 \wedge J_2 \alpha$ by
section \ref{sec:2:3}) won't be solved (without adding others) in any
finite dimensional extension of $G_6$. They are inaudible in the heat
decay anyway. The fact that $E_2$ has no $L^2$ section while plenty of
local ones interprets here as a vanishing theorem similar for instance
of the vanishing of $L^2$ holomorphic functions on $\C$.  Observe
these $J_3$ holomorphic forms were not controlled in the
representation associated to the $T_3$ direction on $Q_7$, but can be
here since we have removed it!

Lastly we complete the study of the other Novikov-Shubin exponents
$\alpha_p$ of $Q_7$. In view of \eqref{eq:2:2:7} the algebraic
pinching theorem \ref{thm:3:3:1} gives $\beta_3 (Q_7) = 2$. By duality
(see discussion around \eqref{eq:G2}) we are left with the study of
$\beta_2 (Q_7) = \beta_4 (Q_7)$. We show that theorem \ref{thm:5:4}
gives $\beta_4 (Q_7)= r_4=1$, that is the system $d_c^{[1]} +
\delta_c$ is C-C elliptic on $E_0^{4,(6)}$.
\begin{proof}
  If $\alpha \in \mathcal{S}_\pi$ belongs to $\ker \pi (d_c^{[1]} +
  \delta_c)$, then by C-C ellipticity of $d_c$, $\beta= \pi (d_c)
  \alpha$ is non trivial if $\alpha$ is.  Moreover $\beta \in
  E_0^{5,(8)} \cap \ker \pi(d_c)$, and this space is Hodge-$*$
  conjugated to $\ker \pi(d^{-*})$ from $\Lambda^{2,-}D^*$ to
  $\Lambda^1 D^*$.  The elliptic symbol (on degenerated
  representations) of this map is well known to be injective. At the
  level of generic representations, one computes easily that
  $\pi(d^{-} d^{-*})$ is an inversible Folland-Stein (or Tanaka)
  operator, giving the result.
\end{proof}
Of course these C-C ellipticity results transplant, like in lemma
\ref{lemma:5:9}, on any elliptic distribution $D^4 \subset TM^7$ with
$M^7$ compact, without integrability condition on the structure.

%\bibliographystyle{plain} 

%\bibliography{tsg}

\noindent
------------------------------------\\
Math\'ematique, B\^at. 425, Universit\'e de Paris-Sud, 91405 Orsay,
France\\
e-mail : \texttt{michel.rumin@math.u-psud.fr}\\
\texttt{http://topo.math.u-psud.fr/\~{}rumin}

\end{document}